\documentclass[12pt,a4paper]{article}
\usepackage {amssymb,latexsym}
\usepackage {amsmath}
\usepackage[french]{babel}

\newtheorem{theorem}{Th\'{e}or\`{e}me}[section]
\newtheorem{prop}{Proposition}[section]
\newtheorem{lemma}{Lemme}[section]
\newtheorem{cor}{Corollaire}[section]

\newtheorem{defi}{D\'{e}finition}[section]
\newtheorem{exer}{Exercice}[subsection]
\newtheorem{prob}{Probl\`{e}me}[section]
\newtheorem{remarque}{Remarque}[section]
\newenvironment{prooof}{
                        \noindent{\bf\small D\'{e}monstration: }\small}
                                       {\hfill {$\mathbf \Box$}\medskip}
\newcommand{\bdefi}{\begin{defi}}
\newcommand{\edefi}{\end{defi}}
\newcommand{\bexer}{\begin{exer}\small\rm}
\newcommand{\eexer}{\end{exer}}
\newcommand{\bsat}{\begin{theorem}}
\newcommand{\esat}{\end{theorem}}
\newcommand{\bprop}{\begin{prop}}
\newcommand{\eprop}{\end{prop}}
\newcommand{\bcor}{\begin{cor}}
\newcommand{\ecor}{\end{cor}}
\newcommand{\blem}{\begin{lemma}}
\newcommand{\elem}{\end{lemma}}
\newcommand{\brem}{\begin{rem}}
\newcommand{\erem}{\end{rem}}
\newcommand{\bbew}{\begin{prooof}}
\newcommand{\ebew}{\end{prooof}}
\newcommand{\bprob}{\begin{prob}}
\newcommand{\eprob}{\end{prob}}
\newcommand{\beq}{\begin{equation}}
\newcommand{\bea}{\begin{eqnarray}}
\newcommand{\eea}{\end{eqnarray}}
\newcommand{\beas}{\begin{eqnarray*}}
\newcommand{\eeas}{\end{eqnarray*}}
\newcommand{\rational}{\mathbb Q}
\newcommand{\real}{\mathbb R}
\newcommand{\korps}{\mathbb K}

\newcommand{\complex}{\mathbb C}
\newcommand{\nat}{\mathbb N}

\newcommand{\lbl}{\label}
\newcommand{\ben}{\begin{enumerate}}
\newcommand{\een}{\end{enumerate}}
\newcommand{\ra}{\rightarrow}
\newcommand{\Cinf}{\mathcal{C}^\infty}

\newcommand{\Ug}{\mathcal{U}\mathfrak{g}}
\newcommand{\Sg}{\mathcal{S}\mathfrak{g}}

\makeatletter
\@addtoreset{equation}{section}

\makeatother

\newcommand{\auteur}[1]{{\sc #1}}

\begin{document}
\thispagestyle{empty}
\begin{center}
%  {\Large \bf D\'{e}formation et rigidit\'{e} des alg\`{e}bres enveloppantes via les
%         structures de Poisson}
{\Large \bf D\'{e}formation par quantification et rigidit\'{e} des alg\`{e}bres
                 enveloppantes}
\vspace{1cm}

     {\bf M.~Bordemann\footnote{{\tt M.Bordemann@uha.fr}},
      A.~Makhlouf\footnote{{\tt N.Makhlouf@uha.fr}},
      T. Petit\footnote{{\tt T.Petit@uha.fr}}}

Laboratoire de Math\'{e}matiques et Applications, \\
Universit\'{e} de Haute Alsace, Mulhouse \\[6mm]
   Novembre 2002
\end{center}
\vspace{7mm}
\small
\begin{center}
 {\bf Abstract}
\end{center}
\begin{center}
\begin{minipage}{10cm}
 We call a finite-dimensional complex
 Lie algebra $\mathfrak{g}$ strongly rigid if its universal enveloping algebra
 $\Ug$ is
 rigid as an associative algebra, i.e. every formal associative deformation
 is equivalent to the trivial deformation. In quantum group theory this
 phenomenon is well-known to be the case for all complex semisimple Lie
 algebras.
 We show that a strongly rigid
 Lie algebra has to be rigid as Lie algebra, and that in addition its
 second scalar cohomology group has to vanish (which excludes nilpotent
 Lie algebras of dimension greater or equal than two). Moreover, using
 Kontsevitch's theory of deformation quantization we show that every
 polynomial deformation of the linear Poisson structure on $\mathfrak{g}^*$
 which induces a nonzero cohomology class of $\mathfrak{g}$ leads to a
 nontrivial deformation of $\mathcal{U}\mathfrak g$. Hence every Poisson
 structure on a vector space which is zero at some point and whose linear
 part is a strongly rigid Lie algebra is therefore formally linearizable
 in the sense of A.~Weinstein. Finally we provide
 examples of rigid Lie algebras which are not strongly rigid, and give a
 classification of all strongly rigid Lie algebras up to dimension 6.
\end{minipage}
\end{center}

\vspace{7mm}

\noindent \textbf{Mots cl\'{e}s}: alg\`{e}bre enveloppante, d\'{e}formation,
quantification, rigidit\'{e}, cohomologie, structure de Poisson et lin\'{e}arisation.


\noindent \textbf{Classification AMS 2000}: 16S30, 16S80, 17Bxx

%\begin{center}
%\textbf{R\'{e}sum\'{e}}
%\end{center}
%
%\begin{center}
%\begin{minipage}{8cm}
%L'objectif de cet article est d'\'{e}tudier les d\'{e}formations et la
%rigidit\'{e} des alg\`{e}bres enveloppantes d'alg\`{e}bres de Lie. Une
%alg\`{e}bre de Lie $\mathfrak{g}$ est dite fortement rigide si son
%alg\`{e}bre enveloppante $\Ug$ est rigide comme
%alg\`{e}bre associative. Notre \'{e}tude montre que la rigidit\'{e} de
%l'alg\`{e}bre de Lie est une condition n\'{e}cessaire pour la rigidit\'{e}
%forte de $\mathfrak{g}$ et qu'il faut se restreindre dans la recherche des
%alg\`{e}bres de Lie fortement rigides aux alg\`{e}bres de Lie dont le
%deuxi\`{e}me groupe de cohomologie scalaire est nul. En utilisant le
%th\'{e}or\`{e}me de formalit\'{e} de Kontsevitsch, nous montrons que toute
%d\'{e}formation polyn\^{o}miale non triviale de la structure de Poisson
%lin\'{e}aire d'une alg\`{e}bre de Lie induit une d\'{e}formation non
%triviale de son alg\`{e}bre enveloppante. Nous donnons des exemples
%d'alg\`{e}bres de Lie rigides qui ne sont pas fortement rigides et la
%classification des alg\`{e}bres de Lie fortement rigides de dimension
%inf\'{e}rieure \`{a} sept.
%\end{minipage}
%\end{center}
\newpage
\normalsize
\thispagestyle{empty}

\tableofcontents

\newpage

\section{Introduction}

Dans la th\'{e}orie des groupes quantiques de Drinfel'd (voir \cite{Dri86}
ou \cite{Kas}), on consid\`{e}re en particulier les
alg\`{e}bres enveloppantes quantiques $\mathcal{U}_h\mathfrak{g}$ d'une
alg\`{e}bre de Lie $(\mathfrak{g},[~,~])$ sur un corps $\korps$
(g\'{e}n\'{e}ralement alg\'{e}briquement clos et de caract\'{e}ristique z\'{e}ro).
Ces alg\`{e}bres sont des d\'{e}formations formelles de l'alg\`{e}bre enveloppante
$\Ug$ de $\mathfrak{g}$ en tant qu'alg\`{e}bre de Hopf. Dans la plupart des
exemples, introduits par Drinfel'D \cite{Dri85} et Jimbo \cite{Jim85},
$\mathfrak{g}$ est une alg\`{e}bre de Lie semi-simple complexe de dimension
finie.
Dans ce cas, la d\'{e}formation formelle de la structure associative de $\Ug$, au sens de
Gerstenhaber \cite{Ge1}, dans
$\mathcal{U}_h\mathfrak{g}$ (ainsi que pour toute autre d\'{e}formation associative)
est une d\'{e}formation triviale, c.-\`{a}-d. isomorphe \`{a} $\Ug$,
(voir par exemple \cite[p.430]{Kas}). En revanche la
comultiplication cocommutative de $\Ug$ se d\'{e}forme en une comultiplication non
cocommutative dans $\mathcal{U}_h\mathfrak{g}$. Ceci signifie que $\Ug$ est
une alg\`{e}bre associative
{\em rigide}, mais non rigide comme alg\`{e}bre de Hopf. Voir \cite{MG} pour les
alg\`{e}bres associatives
complexes rigides de dimension finie. \\
En outre, le th\'{e}or\`{e}me de Poincar\'{e}-Birkhoff-Witt permet de
regarder l'alg\`{e}bre enveloppante $\Ug$ elle-m\^{e}me comme une `d\'{e}formation',
 en g\'{e}n\'{e}ral non commutative, de l'alg\`{e}bre sym\'{e}trique $\Sg$.
Cette derni\`{e}re alg\`{e}bre est isomorphe \`{a} l'alg\`{e}bre des fonctions polyn\^{o}miales
sur l'espace dual $\mathfrak{g}^*$ de $\mathfrak{g}$ et est munie d'une
structure d'alg\`{e}bre de Poisson (provenant de la structure de Poisson
lin\'{e}aire $P_0$ sur $\mathfrak{g}^*$ induite par le crochet de Lie de
$\mathfrak{g}$). La d\'{e}formation $\Sg\leadsto \Ug$
a \'{e}t\'{e} pr\'{e}cis\'{e}e par Gutt \cite{Gu} dans la th\'{e}orie de quantification par
d\'{e}formation de Bayen, Flato, Fr{\o}nsdal, Lichnerowicz et Sternheimer
\cite{BFFLS}. Dans cette th\'{e}orie, on construit des d\'{e}formations
associatives formelles, dites {\em star-produits} $*$, de l'alg\`{e}bre
$\Cinf(M)$ de
toutes les fonctions \`{a} valeurs complexes de classe $\Cinf$ sur une
vari\'{e}t\'{e} de Poisson.
Kontsevitch a r\'{e}examin\'{e} cette d\'{e}formation de $\Sg$ dans le cadre
d'une formule universelle pour un star-produit sur $\real^n$ muni d'une
structure de Poisson arbitraire (voir \cite{Ko}). L'\'{e}quivalence
de la d\'{e}formation de Kontsevitch et celle de Gutt \`{a} \'{e}t\'{e}
d\'{e}montr\'{e}e par Dito
 (\cite{Dit99}).\\
La classification des star-produits \`{a} \'{e}t\'{e} donn\'{e}e par Kontsevitch en termes
de classes de diff\'{e}omorphie formelles des structures de Poisson formelles:
il s'agit des d\'{e}formations formelles des structures de Poisson
elles-m\^{e}mes. En g\'{e}om\'{e}trie diff\'{e}rentielle, ces d\'{e}formations
($\Cinf$ ou formelles) ont \'{e}t\'{e}
implicitement \'{e}tudi\'{e}es dans la th\'{e}orie de lin\'{e}arisation des
structures
de Poisson de Weinstein \cite[pp.~537]{Wei83}: on y a surtout \'{e}tudi\'{e} le cas
o\`{u} la structure de Poisson s'annule en un point et est lin\'{e}arisable dans
le sens qu'elle est isomorphe \`{a} sa partie lin\'{e}aire.
En outre, une formule pour la d\'{e}formation formelle
d'une structure de Poisson de fa\c{c}on que les feuilles symplectiques
restent les m\^{e}mes a \'{e}t\'{e} \'{e}tablie par Lecomte (voir \cite[p.165]{Lec87}).

Toutes ces consid\'{e}rations nous ont amen\'{e}s \`{a} poursuivre deux objectifs
essentiels dans une partie de \cite{Pe} et
dans cet article:\\
D'une part, \'{e}tudier les d\'{e}formations associatives formelles des
alg\`{e}bres enveloppantes, principalement \`{a} l'aide de la d\'{e}formation de la
structure de Poisson lin\'{e}aire sur $\mathfrak{g}^*$ et de la
formule universelle de Kontsevitch. \\
D'autre part, chercher et caract\'{e}riser les alg\`{e}bres de Lie qu'on appelle
{\em fortement rigides}, c.-\`{a}-d. dont
l'alg\`{e}bre enveloppante $\Ug$ est rigide dans le sens qu'elle n'admet que
des d\'{e}formations triviales comme alg\`{e}bre associative.

On a obtenu les r\'{e}sultats principaux suivants:
\begin{enumerate}
 \item A toute d\'{e}formation polyn\^{o}miale formelle non triviale de la structure
   de Poisson lin\'{e}aire provenant d'une alg\`{e}bre de Lie complexe
   $\mathfrak{g}$ on peut associer une d\'{e}formation non triviale de son
   alg\`{e}bre enveloppante $\Ug$, voir le th\'{e}or\`{e}me \ref{Theo4.1}. Pour le
   d\'{e}montrer, on construit d'abord une d\'{e}formation formelle \`{a}
   deux param\`{e}tres de l'alg\`{e}bre sym\'{e}trique $\Sg$ d'apr\`{e}s Kontsevitch
   et on montre sa convergence par rapport au premier param\`{e}tre.
 \item Nous avons montr\'{e}, ind\'{e}pendamment de la
   th\'{e}orie de Kontsevitch, les deux cas
   particuliers importants du th\'{e}or\`{e}me \ref{Theo4.1} suivants:
   \begin{itemize}
    \item Si l'alg\`{e}bre de Lie $\mathfrak{g}$ n'est pas rigide,
        alors $\Ug$ n'est pas
        rigide non plus, voir th\'{e}or\`{e}me \ref{Theo1.3}. Ceci justifie
        la notion `fortement rigide' et nous permet d'utiliser les
        r\'{e}sultats sur la rigidit\'{e} des alg\`{e}bres de Lie obtenus par
        R.~Carles, Y.~Diakit\'{e}, M.~Goze et J.M.~Ancochea-Berm\'{u}dez
        dans (\cite{Ca-Di},\cite{Ca1}, \cite{GA}, \cite{AB-Go}).
    \item Si le deuxi\`{e}me groupe de cohomologie scalaire,
       $\mathbf{H}_{CE}^2(\mathfrak{g},\korps)$, de l'alg\`{e}bre de Lie
       $\mathfrak{g}$ ne s'annule pas,
       alors $\Ug$ n'est pas rigide, voir th\'{e}or\`{e}me \ref{Theo2.1}.
       Ce crit\`{e}re montre par exemple que les alg\`{e}bres de Lie nilpotentes de
       dimension $\geq 2$ ne sont pas fortement rigides.
   \end{itemize}
 \item En plus des exemples d\'{e}j\`{a} connus comme la classe des alg\`{e}bres de Lie
    semi-simples et l'alg\`{e}bre non ab\'{e}lienne $\mathfrak{r}_{2}$
    de dimension $2$, on a \'{e}tabli la forte rigidit\'{e}
    des alg\`{e}bres affines $\mathfrak{gl}(m,\korps)\times \korps^m$
     (voir le th\'{e}or\`{e}me \ref{TInfFortRigGlmKm}). Une preuve
      plus g\'{e}om\'{e}trique \`{a} \'{e}t\'{e} ind\'{e}pendamment obtenue par Dufour et Zung
    \cite{DNTZ02} dans le cadre des structures de Poisson lin\'{e}arisables.
     Pour tous ces exemples, nous avons v\'{e}rifi\'{e} la
    nullit\'{e} du deuxi\`{e}me groupe de cohomologie
    $\mathbf{H}_{CE}^2(\mathfrak{g},\Sg)$ qui est toujours isomorphe au
    deuxi\`{e}me groupe de cohomologie de Hoch\-schild
    $\mathbf{H}_{H}^2(\Ug,\Ug)$ d'apr\`{e}s un th\'{e}or\`{e}me de H.~Cartan et
    S.~Eilenberg.
 \item Notre principal th\'{e}or\`{e}me  \ref{Theo4.1} nous permet d'exclure
    quelques alg\`{e}bres de Lie rigides de la liste des fortement rigides,
    par exemple les sommes semi-directes
    $\mathfrak{t}_{1}\times \mathfrak{n}_{5,6}$ (alg\`{e}bre de Lie r\'{e}soluble) et
    $\mathfrak{sl}(2,\complex)\times \mathfrak{n}_3$ (o\`{u} $\mathfrak{n}_3$
    est l'alg\`{e}bre de Heisenberg). Ainsi, on classifie
    les alg\`{e}bres de Lie fortement rigides jusqu'en dimension $6$:
    \begin{center}
        $\{0\},~~\korps,~~\mathfrak{r}_{2}, ~~\mathfrak{sl}(2,\korps),~~
           \mathfrak{gl}(2,\korps),~~
            \mathfrak{sl}(2,\korps)\times \mathfrak{r}_{2},~~$\\
        $\mathfrak{sl}(2,\korps)\times \mathfrak{sl}(2,\korps),~~
              \mathfrak{gl}(2,\korps)\times \korps^{2}$.
    \end{center}

\item Du theor\`{e}me \ref{Theo4.1}, on en d\'{e}duit \'{e}galement que toute
    structure de Poisson dont la partie lin\'{e}aire correspond \`{a} une
    alg\`{e}bre fortement rigide est formellement lin\'{e}arisable.
\end{enumerate}

Dans cet article, on consacre le deuxi\`{e}me paragraphe aux rappels de quelques g\'{e}n\'{e}ralit\'{e}s
concernant
les alg\`{e}bres enveloppantes, la th\'{e}orie des d\'{e}formations, la cohomologie,
la rigidit\'{e} et les structures de Poisson.
Dans le paragraphe 3, on d\'{e}finit les alg\`{e}bres de Lie fortement rigides,
et on trouve les exemples mentionn\'{e}s ci-dessus par un calcul
cohomologique.
%on dira qu'une alg\`{e}bre de Lie est fortement rigide si son alg\`{e}bre enveloppante
%est rigide comme  alg\`{e}bre associtive. C'est \`{a} dire que toute d\'{e}formation est
%isomorphe \`{a} une d\'{e}formation triviale. On donne, ensuite, comme exemples les
%alg\`{e}bres de
% Lie semi-simples (voir par exemple \cite[p.430]{Kas}), l'alg\`{e}bre r\'{e}soluble
% affine de dimension 2 et
%la somme semi-directe  $\mathfrak{gl}(n,\korps)\times \korps^n$. On montre
%que pour tous
%ces exemples le deuxi\`{e}me groupe de cohomologie de Hochschild s'annule.
Le paragraphe 4 est consacr\'{e} aux propri\'{e}t\'{e}s des alg\`{e}bres de Lie fortement
rigides, notamment \`{a} la d\'{e}\-monstration du th\'{e}or\`{e}me \ref{Theo1.3}
et \ref{Theo2.1} mentionn\'{e}s ci-dessus.
% On \'{e}tablit qu'une alg\`{e}bre enveloppante rigide
%provient toujours d'une alg\`{e}bre de Lie rigide. Ceci restreint la recherche des
%alg\`{e}bres
% de Lie fortement rigides aux seules
%alg\`{e}bres de Lie rigides. Ensuite, on montre que si le deuxi\`{e}me groupe de
%cohomologie
%scalaire d'une alg\`{e}bre de Lie est non nul, alors elle n'est pas fortement
%rigide. Ce dernier r\'{e}sultat permet d'\'{e}tablir une liste d'alg\`{e}bres de Lie
%rigides qui ne sont
%pas fortement rigides.
Dans le paragraphe 5, on rappelle la th\'{e}orie de la quantification par
d\'{e}formation et la formule universelle de
Kontsevitsch pour une d\'{e}formation de $\Cinf(M)$ pour toute vari\'{e}t\'{e} de Poisson
$M$. Apr\`{e}s avoir rappel\'{e} la construction de Gutt et de Dito, on d\'{e}montre
notre th\'{e}or\`{e}me principal \ref{Theo4.1} mentionn\'{e} ci-dessus.
%d\'{e}formation polyn\^{o}miale non triviale de la structure de Poisson
%lin\'{e}aire d'une alg\`{e}bre de Lie (\ref{Eq(2)}) induit une d\'{e}formation non
%triviale de son alg\`{e}bre enveloppante. Ce
%th\'{e}or\`{e}me sera utilis\'{e} pour montrer qu'il existe des
%alg\`{e}bres de Lie rigides dont le deuxi\`{e}me groupe de cohomologie
%scalaire est nul et qui ne sont pas fortement rigides.
Dans le paragraphe 6,
on donne
la classification des alg\`{e}bres de Lie fortement rigides de
dimension inf\'{e}rieure ou \'{e}gale \`{a} six.
On en deduit \'{e}galement quelques
remarques a propos du probl\`{e}me de lin\'{e}arisation des structures de Poisson
en paragraphe 7.

\subsubsection* {Remerciments}

\noindent Nous remercions Aleksei Bolsinov, Roger Carles,
Jean-Paul Dufour, Michel Goze et Freddy Van Oystaeyen pour les
fructueuses discussions et propositions.

\section{G\'{e}n\'{e}ralit\'{e}s}
\begin{enumerate}

\item Soient $K$ un anneau commutatif unitaire et $\mathfrak{g}$
une alg\`{e}bre de Lie
sur $K$. Rappelons qu'une $\mathfrak{g}$-repr\'{e}sentation (\`{a} gauche)
de $\mathfrak{g}$ est un $K$-module $\mathcal{M}$ et un
$K$-homomorphisme
\beq
\mathfrak{g}\otimes_K \mathcal{M} \mapsto \mathcal{M}:
x\otimes a \mapsto xa
\end{equation}
tel que $x(ya)-y(xa)=[x,y]a$.
 A toute alg\`{e}bre de Lie $\mathfrak{g}$,
on associe une $K$-alg\`{e}bre associative $\Ug$
telle que toute
 $\mathfrak{g}$-repr\'{e}sentation
(\`{a} gauche) peut \^{e}tre vue comme une $\Ug$-repr\'{e}sentation
(\`{a} gauche)
et vice-versa. L'alg\`{e}bre $\Ug$ est construite de la
mani\`{e}re suivante:

Soit $T\mathfrak{g}$ l'alg\`{e}bre tensorielle du $K$-module
$\mathfrak{g}$ ,
$T\mathfrak{g}=T^{0}\oplus T^{1}\oplus \cdots \oplus T^{n}\oplus \cdots$ o\`{u}
$T^{n}=\mathfrak{g}\otimes_K \mathfrak{g}\otimes_K \cdot \cdot \cdot
\otimes_K
\mathfrak{g}$ ($n$ fois),
en particulier $T^{0}=K\hspace{0.5mm}1$ et $T^{1}=\mathfrak{g}$.
Le produit
dans $T\mathfrak{g}$ \'{e}tant la multiplication tensorielle.
Toute application
$K$-lin\'{e}aire $\mathfrak{g}\otimes_K \mathcal{M}\rightarrow \mathcal{M}$
admet une extension
unique en un morphisme de $K$-alg\`{e}bres $T\mathfrak{g}\rightarrow
\mathbf{Hom}_K(\mathcal{M},\mathcal{M})$. Pour que la restriction de
ce morphisme \`{a} $\mathfrak{g}$ d\'{e}finisse une $\mathfrak{g}$-repr\'{e}sentation
sur $\mathcal{M}$,
%le $\mathfrak{g}$-module (?)
il est n\'{e}cessaire
et suffisant que les \'{e}l\'{e}ments de $T\mathfrak{g}$ de la forme
$x\otimes y-y\otimes x-[x,y]$
o\`{u} $x,y\in \mathfrak{g}$ annulent $\mathcal{M}$. Par cons\'{e}quent,
on introduit l'id\'{e}al
bilat\`{e}re $\mathcal{I}$ engendr\'{e} par les \'{e}l\'{e}ments $x\otimes y-y\otimes x-[x,y]$
o\`{u} $x,y\in \mathfrak{g}$, et on d\'{e}finit l'alg\`{e}bre enveloppante de
$\mathfrak{g}$
comme le quotient $T\mathfrak{g}/\mathcal{I}$. On identifie ainsi les
$\mathfrak{g}$-repr\'{e}sentations
et les $\Ug$-modules. Rappelons que tout $\Ug$-bimodule
$\mathcal{M}$ est un $\mathfrak{g}$-module par $(x,m)\rightarrow xm-mx$,
not\'{e} $\mathcal{M}_a$.

Supposons que l'alg\`{e}bre de Lie $\mathfrak{g}$ soit un $K$-module libre, soit
$\{x_{i}\}$
une base fix\'{e}e de $\mathfrak{g}$. Soit $y_{i}$ l'image de $x_{i}$
par le $K$-homomorphisme
$\mathfrak{g}\rightarrow \Ug$. On pose
$y_{I}=y_{i_{1}}\cdots y_{i_{p}}$ avec $I$ une suite finie
d'indices $1\leq i_{1}\leq \ldots\leq i_{p}$ et $y_{I}=1$ si $I=\emptyset $.
D'apr\`{e}s le th\'{e}or\`{e}me de Poincar\'{e}-Birkhoff-Witt, l'alg\`{e}bre enveloppante
$\Ug$ est engendr\'{e}e par les \'{e}l\'{e}ments $y_{I}$
correspondant
\`{a} des suites croissantes $I$ (voir par exemple \cite[pp.~271]{Ca-Ei}).

Soit $K=\korps$ un corps. On d\'{e}signe par $\mathcal{S}V$ l'alg\`{e}bre sym\'{e}trique
sur le $\korps$-module $V$.
Si $\rational \in \korps$, il existe une bijection canonique
$\omega:\Sg\ra\Ug$ donn\'{e}e par
\beq\lbl{EqDefSymmSgUg}
  \omega(x_1\cdots x_n):=\frac{1}{n!}\sum_{\sigma\in S_n}
                             x_{\sigma(1)}\cdots x_{\sigma(n)}
\end{equation}
(o\`{u} $x_1,\cdots,x_n\in\mathfrak{g}$ et $S_n$ est le groupe des
permutations de $\{1,\ldots,n\}$)
qui est un
 isomorphisme de $\mathfrak{g}$-modules entre $\Sg_{a}$ et
 $\Ug_a$ (voir par exemple \cite[pp.78-79]{Di2}).

\item
Les premiers travaux sur les d\'{e}formations des alg\`{e}bres remontent aux
ann\'{e}es 60. L'article de M. Gerstenhaber \cite{Ge1} fournit l'outil de ces
d\'{e}\-for\-ma\-tions sur la base des s\'{e}ries formelles qu'on rappelle
ci-dessous.\\
D\'{e}sormais $\korps$ d\'{e}signe un corps alg\'{e}briquement clos de caract\'{e}ristique
nulle. Soit  $\korps[[t]]$ l'anneau des s\'{e}ries formelles \`{a}
coefficients
dans $\korps$. Pour chaque $\korps$-espace vectoriel $E$ on note
$E[[t]]$ le $\korps[[t]]$-module de toutes les s\'{e}ries formelles \`{a} coefficients
dans $E$. Soit $(\mathcal{A},\mu_{0})$ une $\korps$-alg\`{e}bre associative
(resp. de Lie).
Alors $(\mathcal{A}[[t]],\mu_{0})$ est une $\korps[[t]]$-alg\`{e}bre
associative (resp. de Lie).

 \bdefi \lbl{Def1.1}

 \begin{enumerate}
   \item On appelle {\em d\'{e}formation associative} (resp. {\em d'une alg\`{e}bre
  de Lie}) {\em formelle} de $\mathcal{A}$ une
   $\korps[[t]]$-alg\`{e}bre associative (resp. de Lie)
   $(\mathcal{A}[[t]],\mu_{t})$
  avec
  \[
    \mu _{t}=\mu_{0}+t\mu_{1}+t^{2}\mu_{2}+\cdots +t^{n}\mu_{n}+\cdots,
  \]
  o\`{u} $\mu_{n}\in \mathbf{Hom}_{\korps}(\mathcal{A}\otimes_{\korps}
  \mathcal{A},\mathcal{A}).$ (resp. $\mu_{n}\in \mathbf{Hom}_{\korps}
  (\mathcal{A}\wedge_{\korps}\mathcal{A},\mathcal{A}).$

   \item Deux d\'{e}formations $(\mathcal{A}[[t]],\mu_t)$ et
  $(\mathcal{A}[[t]],{\mu}_t^\prime)$
  sont dites {\em \'{e}quivalentes} s'il existe un isomorphisme formel
  \[
     \varphi_{t}=\varphi _{0}+\varphi _{1}t+\cdots+\varphi _{n}t^{n}+\cdots,
  \]
  avec $\varphi_{0}=Id_{\mathcal{A}}$ (l'application identit\'{e} de
  $\mathcal{A}$)
  et $\varphi_{n}\in \mathbf{Hom}_{\korps}(\mathcal{A},\mathcal{A})$
  tel que
  \[
       \mu_{t}^{\prime}(a,b)=
           \varphi _{t}^{-1}(\mu_{t}(\varphi_{t}(a),\varphi_{t}(b))
                  ~~~\forall a,b\in \mathcal{A.}
  \]

   \item Une d\'{e}formation de $\mathcal{A}$ est dite {\em triviale} si elle
  est \'{e}quivalente \`{a} $(\mathcal{A}[[t]],\mu_0)$.

   \item Une alg\`{e}bre associative (resp. de Lie) $\mathcal{A}$ est dite
   {\em rigide} si toute d\'{e}formation de $\mathcal{A}$ est triviale.
  \end{enumerate}
  \edefi

  \item La th\'{e}orie des d\'{e}formations est intimement li\'{e}e \`{a} la cohomologie de
  Hochschild dans le cas des alg\`{e}bres associatives et \`{a} la cohomologie de
  Chevalley-Eilenberg dans le cas des alg\`{e}bres de Lie. On d\'{e}signe par
  $\mathbf{H}_{H}^{n}(\mathcal{A},\mathcal{M})$ le $n$-i\`{e}me groupe de
  cohomologie de Hochschild d'une alg\`{e}bre associative $\mathcal{A}$ \`{a} valeurs
  dans un bimodule $\mathcal{M}$ et par
  $\mathbf{H}_{CE}^{n}(\mathfrak{g},\mathcal{M})$ le $n$-i\`{e}me groupe de
  cohomologie de Chevalley-Eilenberg d'une alg\`{e}bre de Lie $\mathcal{A}$ \`{a}
  valeurs dans un $\mathfrak{g}$-module $\mathcal{M}$.
   Le deuxi\`{e}me groupe de cohomologie de Hochschild d'une alg\`{e}bre associative
   (resp. cohomologie de Chevalley-Eilenberg d'une alg\`{e}bre de Lie) \`{a} valeurs
   dans l'alg\`{e}bre correspond \`{a} l'espace des d\'{e}formations infinit\'{e}simales.
   Pour notre \'{e}tude, le fait suivant sera tr\`{e}s important (voir \cite{Ge1},
  \cite{Ge2} et \cite{Ni-Ri}):
  \bprop \lbl{PHDeuxNulRig}
   Soit $\mathcal{A}$ (resp. $\mathfrak{g}$) une alg\`{e}bre associative
   (resp. une alg\`{e}bre de Lie).
   \ben
   \item Si le deuxi\`{e}me groupe de cohomologie de Hochschild,
   $\mathbf{H}_H^2(\mathcal{A},\mathcal{A})$,
   (resp. de Chevalley-Eilenberg,
   $\mathbf{H}_{CE}^2(\mathfrak{g},\mathfrak{g})$,)
   s'annule, alors l'alg\`{e}bre $\mathcal{A}$ (resp. $\mathfrak{g}$) est rigide.
   \item D'un autre cot\'{e}, si $\mu_0+t\mu_1+\ldots$ est une d\'{e}formation de
   $\mathcal{A}$ (resp. de $\mathfrak{g}$) et si  $\mathcal{A}$
   (resp. $\mathfrak{g}$)
   est rigide, alors le $2$-cocycle $\mu_1$ est toujours un $2$-cobord de
   la cohomologie
   de Hochschild (resp. de Chevalley-Eilenberg).
   \een
   \eprop
   (voir \cite[p.430]{Kas} pour une d\'{e}monstration). La r\'{e}ciproque de
   l'assertion (a) est fausse, il existe des exemples d'alg\`{e}bres de Lie rigides
   dont le deuxi\`{e}me groupe de cohomologie est non nul (voir \cite{Ri}. Le
   probl\`{e}me est encore ouvert en caract\'{e}ristique $0$ pour les alg\`{e}bres
   associatives (voir \cite{Ge-Ds} pour un exemple en caract\'{e}ristique $p$)

  Le troisi\`{e}me groupe de cohomologie quand \`{a} lui correspond \`{a} l'espace des
  obstructions pour l'extension d'ordre $n$ \`{a} une d\'{e}formation d'ordre $n+1$
  (\cite{Ge1}, \cite{Ge2} et \cite{Ni-Ri}).

\item
Le th\'{e}or\`{e}me classique suivant d\^{u} \`{a} H.~Cartan et S.~Eilenberg
(\cite[pp.277]{Ca-Ei}) fait le lien entre la cohomologie de Hochschild de
l'alg\`{e}bre enveloppante \`{a} valeurs dans un $\Ug$-bimodule
$\mathcal{M}$ (en particulier $\mathcal{M}=
\Ug$) et la cohomologie de Chevalley-Eilenberg de l'alg\`{e}bre
de Lie \`{a} valeurs dans le m\^{e}me module.

\bsat \lbl{TheoC-E}
Soit $\mathfrak{g}$ une alg\`{e}bre de Lie de dimension finie sur $\korps$.  Alors
\[
   \mathbf{H}_{H}^{n}(\Ug,\mathcal{M})
       \simeq \mathbf{H}_{CE}^{n}(\mathfrak{g},\mathcal{M}_{a})
            ~~~\forall n\in \nat
\]
En particulier, si $\mathbb{Q} \subset \korps$ on a pour tout $n\in \nat$
\[
   \mathbf{H}_{H}^{n}(\Ug,\Ug)
       \simeq \mathbf{H}_{CE}^{n}(\mathfrak{g},\Ug_{a})
       \simeq \mathbf{H}_{CE}^{n}(\mathfrak{g},\Sg_{a})
       \simeq \bigoplus_{k=0}^\infty
           \mathbf{H}_{CE}^{n}(\mathfrak{g},\mathcal{S}^k\mathfrak{g}_{a})
\]
\esat

%\vspace{0.03in}\textbf{Th\'{e}or\`{e}me (H.S). }
\item Le th\'{e}or\`{e}me de Hochschild et Serre (\cite{Ho-Se})
suivant fournit une factorisation des
groupes de cohomologie de Chevalley-Eilenberg :
\bsat  \lbl{TheoH.S}
Soient $\mathfrak{g}$ une alg\`{e}bre de Lie de dimension finie, $\mathcal{M}$
un $\mathfrak{g}$-module de dimension finie sur $\korps$,
$\mathfrak{n}$ un id\'{e}al de $\mathfrak{g}$ et $\mathfrak{b}$ une
sous-alg\`{e}bre suppl\'{e}mentaire de $\mathfrak{n}$,
r\'{e}ductive dans $\mathfrak{g}$, telle que le $\mathfrak{b}$-module
induit sur $\mathcal{M}$ soit semi-simple, alors
pour tout entier positif $p$ on a
\[
   \mathbf{H}_{CE}^{p}(\mathfrak{g},\mathcal{M})\simeq
       \sum_{i+j=p}\mathbf{H}_{CE}^{i}(\mathfrak{b},\korps)
              \otimes
              \mathbf{H}_{CE}^{j}(\mathfrak{n},\mathcal{M})^{\mathfrak{b}}.
\]
o\`{u} $\mathbf{H}_{CE}^{j}(\mathfrak{n},\mathcal{M})^{\mathfrak{b}}$ d\'{e}signe le
sous-espace de tous les \'{e}l\'{e}ments invariants par $\mathfrak{b}$.
\esat

\item

Une {\em alg\`{e}bre de Poisson} est une
alg\`{e}bre associative commutative $\mathcal{A}$ sur $\korps$ munie
 d'une application bilin\'{e}aire $\{~,~\}:\mathcal{A}\times \mathcal{A}\rightarrow
 \mathcal{A}$
v\'{e}rifiant pour tout $f,g,h\in \mathcal{A}$
\bea
    \{f,g\} & = & -\{g,f\} \lbl{Pois1} \\
0 & = & \{f,\{g,h\}\}+\{g,\{h,f\}\}+\{h,\{f,g\}\} \nonumber \\
                & &   \qquad \textit{(identit\'{e} de Jacobi)} \lbl{Pois2}
                             \\
    \{h,fg\}& = & \{h,f\}g+f\{h,g\} \nonumber \\
                & & \qquad \textit{(identit\'{e} de Leibniz)}
                      \lbl{Pois3}
\eea
On note par $\left( \mathcal{A},\cdot ,\{~,~\}\right)$ une telle alg\`{e}bre.
Une vari\'{e}t\'{e} $M$ est dite {\em vari\'{e}t\'{e} de Poisson} si l'alg\`{e}bre des fonctions
$\Cinf\left(M\right)$ est munie d'une structure de Poisson.

\item  Soit $M$ une vari\'{e}t\'{e} de Poisson.
Pour tout $f\in \Cinf\left(M\right)$, l'application
lin\'{e}aire $\{f,~\}:\Cinf\left(M\right) \rightarrow
\Cinf\left(M\right)$ d\'{e}finit une d\'{e}rivation. Soit
$X_{f}$ le champ de vecteurs associ\'{e} \`{a} $f$, le crochet de Poisson
s'\'{e}crit alors
\[
 \{f,g\}=X_{f}g=-X_{g}f=dg(X_{f})=-df(X_{g})
\]

\vspace{0.05in}

Il s'ensuit que le crochet de Poisson est d\'{e}termin\'{e} par une
application bilin\'{e}aire antisym\'{e}trique sur le fibr\'{e} cotangent
 $T^{*}M$ de $M$, c'est \`{a} dire un champ de
 tenseurs antisym\'{e}triques $P\in \Gamma^\infty(M,\Lambda^{2}TM)$
 (avec $TM$ le fibr\'{e} tangent de $M$) tel que
\beq
   \{f,g\}=P(df,dg)=\sum_{i,j}P^{ij}\partial_{i}f\partial_{j}g
\end{equation}
avec $\partial_{i}$ $(1\leq i\leq n)$ la d\'{e}rivation
$\frac{\partial}{\partial x^{i}}$ o\`{u}
$\left(x^1,\ldots,x^n\right)$ sont des coordonn\'{e}es locales sur
$M$ et $P^{ij}$ sont des fonctions
locales $\Cinf$ d\'{e}finies par $P^{ij}:=\{x^i,x^j\}$.
Le champ de tenseurs $P$ est appel\'{e} {\em bivecteur de Poisson}.

Toute structure de Poisson sur $M$ est \'{e}quivalente \`{a}
la donn\'{e}e d'un bivecteur $P\in \Gamma^\infty(M,\Lambda^{2}TM)$
v\'{e}rifiant, dans toute carte,
\beq
  \sum_{h=1}^n \big( P^{ih}\partial_{h}P^{jk}
                     +P^{jh}\partial_{h}P^{ki}
                     +P^{kh}\partial_{h}P^{ij}\big)=0.
\end{equation}

\item
Soient $M$ une vari\'{e}t\'{e} diff\'{e}rentiable,
$\mathcal{T}^{i}M:=\Gamma^\infty(M,\Lambda^{i}TM)$
l'espace des champs de tenseurs antisym\'{e}triques de degr\'{e} $i\in\nat$,
$\mathcal{T}^{0}M:=\Cinf (M)$ et
$\mathcal{T}M=(\oplus_{n\geq 0}\mathcal{T}^{i}M,\wedge)$
l'alg\`{e}bre des multivecteurs sur $M$.

Pour tout $X,Y_{1},\cdots ,Y_{p} \in \mathcal{T}^{1}M$ l'op\'{e}rateur
d\'{e}riv\'{e}e de Lie associ\'{e} est d\'{e}fini par $\;\;$\ $\;\;\;$
\beq
[X,Y_{ 1}\wedge \cdots \wedge Y_{p}]=\sum_{i=1}^{p}Y_{1}\wedge \cdots \wedge
[X,Y_{i}] \wedge \cdots \wedge Y_{p}
\end{equation}
On d\'{e}finit pour $X_{1},\cdots ,X_{p}\in \mathcal{T}^{1}M$ et $Q\in
\mathcal{T}M$ une extension de la d\'{e}riv\'{e}e de Lie,
appel\'{e}e {\em crochet de Schouten}, d\'{e}finie par:
\beq
  [X_{1}\wedge \cdots \wedge X_{p},Q]_{s}
   =\sum_{i=1}^{p}\left(-1\right)^{i+1}X_{1}\wedge \cdots \wedge
\widehat{X_{i}}\wedge \cdots X_{p}\wedge [X_{i},Q]
\end{equation}
avec $\widehat{X_{i}}$ signifiant que $X_{i}$ est omis.

\vspace{0.05in}

Le crochet de Schouten
$[~,~]_{s}:\mathcal{T}^{p}M\times \mathcal{T}^{q}M
\rightarrow \mathcal{T}^{p+q-1}M$ v\'{e}rifie les
propri\'{e}t\'{e}s suivantes :
\bea
~  [P,Q]_{s} & = & -(-1)^{(p-1)(q-1)}[Q,P]_{s} \lbl{Schou1} \\
~ [P,Q\wedge R]_{s} & = & [P,Q]_{s}\wedge R+(-1)^{pq+q}Q\wedge [P,R]_{s}
                              \lbl{Schou2} \\
  0 & = & (-1)^{(p-1)(r-1)}[P,[Q,R]_{s}]_{s}+(-1)^{(q-1)(p-1)}[Q,[R,P]_{s}]_{s}
                     \nonumber \\
    & & ~~~~~~+(-1)^{(r-1)(q-1)}[R,[P,Q]_{s}]_{s}  \lbl{Schou3}
\eea

L'alg\`{e}bre des multivecteurs $\mathcal{T}M$ munie du crochet de
Schouten est une super-alg\`{e}bre de Lie. En particulier,
on a les formules suivantes en coordonn\'{e}es locales:
 \ben
 %$(1)$
 \item
 Soient
$P=\frac{1}{2}\sum_{i,j}P^{ij}\partial_{i}\wedge \partial_{j}$
($\in \mathcal{T}^{2}M$) et
$Q=\frac{1}{2}\sum_{i,j}Q^{ij}\partial_{i}\wedge \partial_{j}$
($\in \mathcal{T}^{2}M$). Il vient:
\begin{eqnarray}
[P,Q]_{s}   & = &  \sum_{i,j,k,h}
                     \big( P^{ih}\partial_{h}Q^{jk}
                        +P^{kh}\partial_{h}Q^{ij}
                        +P^{jh}\partial _{h}Q^{ki} \nonumber \\
            &   &        ~~~~~~~~+Q^{ih}\partial_{h}P^{jk}
                         +Q^{kh}\partial_{h}P^{ij}
                         +Q^{jh}\partial_{h}P^{ki}\big)
                         \partial_{i}\wedge \partial_{j}\wedge \partial_{k}
                         \nonumber \\
            &   &
\end{eqnarray}
%$(2)$
 \item
 Soient $P=\frac{1}{2}\sum_{i,j}P^{ij}\partial_{i}\wedge \partial_{j}$
($\in \mathcal{T}^{2}M$) un champ de bivecteurs et
$A=\sum_{i}A^{i}\partial _{i}$
($\in \mathcal{T}^{1}M$) un champ de vecteurs. Il vient:
\beq
    [P,A]_{s}=\sum_{i,j,k}\big(
         P^{ik}\partial_{k}A^{j}
         +P^{kj}\partial_{k}A^{i}
         -P^{ij}\partial_{k}A^{k}\big)
         \partial_{i}\wedge \partial _{j}
\end{equation}
\een

Un bivecteur $P\in \mathcal{T}^{2}M$ d\'{e}finit une structure de
Poisson si et seulement si le crochet de Schouten $[P,P]_{s}=0$. Il
s'ensuit que l'op\'{e}rateur $\delta_{P}:=[P,~]_{s}$ est de carr\'{e} nul et d\'{e}termine
la structure d'un complexe sur $\mathcal{T}M$ dont la cohomologie est
appel\'{e}e {\em cohomologie de Poisson}.

Soit $P_t:=\sum_{a=0}^\infty t^a P_a$ une s\'{e}rie formelle de champs de
bivecteurs, c.-\`{a}-d. $P_a\in\mathcal{T}^{2}M$ quel que soit
$a\in\nat$. $P_t$ est dite {\em structure de Poisson formelle sur $M$}
lorsque $[P_t,P_t]_s=0$ ce qui est \'{e}quivalent \`{a}
\beq \lbl{EqDefPoisDefo}
    \sum_{\stackrel{b,c\geq 0}{b+c=a}}[P_b,P_c]_{s} = 0
     ~~~~~~~\forall a\in\nat
\end{equation}
En particulier, $[P_0,P_0]_s=0$, alors $P_0$ est une structure de Poisson
sur $M$. La s\'{e}rie $P_t$ est dite {\em d\'{e}formation formelle de la
structure de Poisson $P_0$}. Pour $a=1$, l'\'{e}quation (\ref{EqDefPoisDefo})
donne $0=[P_0,P_1]_s$, donc le champ de bivecteurs $P_1$ est toujours un
cocycle de la cohomologie de Poisson.

\item

Soit $(\mathfrak{g},[~,~])$ une alg\`{e}bre de Lie de dimension finie $n$
sur un corps $\korps$ et $\mathfrak{g}^{*}$ son dual
alg\'{e}brique. L'alg\`{e}bre sym\'{e}trique $\Sg$ s'identifie de mani\`{e}re naturelle
\`{a} l'alg\`{e}bre des fonctions polyn\^{o}miales sur $\mathfrak{g}^{*}$.
La structure de $\mathfrak{g}$ permet de d\'{e}finir un crochet de
Poisson lin\'{e}aire sur $\mathfrak{g}^{*}$:
\[
   \{f,g\}(x):=x([df(x),dg(x)])
\]
avec $f,g \in \Sg$ et $x \in \mathfrak{g}^{*}$.
Soient $(e_{i})_{i=1,\ldots, n}$ une base de
$\mathfrak{g}$, $(e^{i})_{i=1,\ldots,n}$ sa base duale et
$x =\sum_{i=1}^{n}x_{i}e^{i}\in \mathfrak{g}^{*}$,
$\{f,g\}=P_{0}(df,dg)$ avec
$P_{0}$ le bivecteur d\'{e}fini par
\beq\lbl{Eq(2)}
        P_{0}=
          \frac{1}{2}\sum_{i,j}
             P_{0}^{ij}\partial_{i}\wedge \partial_{j}
             \mathrm{~~~o\grave{u}~~}
                P_{0}^{ij}(x)=\sum_{k}C_{ij}^{k}x_{k}
\end{equation}
avec $C_{ij}^{k}$ les constantes de structure de $\mathfrak{g}$.

Il s'ensuit que $(\Sg,\cdot,\{~,~\})$ est une alg\`{e}bre de Poisson.

\item Toute cocha\^{\i}ne de Chevalley-Eilenberg dans
$\mathbf{Hom}_\korps(\Lambda^{k}\mathfrak{g},\Sg)$ s'identifie de
fa\c{c}on naturelle \`{a} un champ de
multivecteurs dans  $\mathcal{T}^{k}(\mathfrak{g}^{*})$.
Les deux formules suivantes qui se d\'{e}duisent par un calcul direct
caract\'{e}risent la liaison entre l'op\'{e}rateur cobord $\delta_{CE}$ de la
cohomologie de
Chevalley-Eilenberg et celui de la structure de Poisson lin\'{e}aire.
\blem \lbl{lm}
Soit $\phi \in \mathbf{Hom}_\korps(\Lambda^{k}\mathfrak{g},\Sg)
   \cong \Lambda^{k}\mathfrak{g}^{*}\otimes \Sg
   \subset \mathcal{T}^{k}(\mathfrak{g}^{*})$ et soit
   $i(e_{i}):\Lambda^{\bullet}\mathfrak{g}\ra\Lambda^{\bullet-1}\mathfrak{g}$
   le produit int\'{e}rieur.
\bea
 \delta_{CE}\phi (x) & = &
        \sum_{i,j,k}C_{ik}^{j}x_{j}e^{i}\wedge
                     \frac{\partial}{\partial {x_{k}}}\phi(x)
               -\frac{1}{2}\sum_{i,k,l}
               C^{i}_{kl}e^{k}\wedge e^{l}\wedge i(e_{i})(\phi(x))
                   \nonumber \\
                     &   & \lbl{EqCheEiFormule} \\
 \delta_{CE}\phi     & = & [P_{0},\phi]_s \lbl{EqCheEiSchouten}
\eea
\elem

\end{enumerate}

\section{Alg\`{e}bres de Lie fortement rigides}

\noindent Dans ce paragraphe on introduit la notion de forte rigidit\'{e} avec
quelques exemples.

\bdefi \lbl{Def1.2}
Une alg\`{e}bre de Lie $\mathfrak{g}$ est dite {\em fortement
rigide} si son alg\`{e}bre enveloppante $\Ug$ est rigide comme
alg\`{e}bre associative.\\
Une alg\`{e}bre de Lie  $\mathfrak{g}$ est dite {\em infinit\'{e}simalement fortement
rigide} si le deuxi\`{e}me groupe de cohomologie
$\mathbf{H}^2_{CE}(\mathfrak{g},\Sg)$ s'annule.
\edefi
La proposition suivante est une cons\'{e}quence imm\'{e}diate de la proposition
\ref{PHDeuxNulRig} et du th\'{e}or\`{e}me de
Cartan-Eilenberg (\ref{TheoC-E}):
\bprop
 Soit $\mathfrak{g}$ une alg\`{e}bre de Lie de dimension finie. Si
 $\mathfrak{g}$ est infinit\'{e}simalement fortement rigide, alors elle est
 fortement rigide.
\eprop
Tous les exemples qui suivent correspondent \`{a} des alg\`{e}bre de Lie
in\-fin\-it\'{e}\-si\-male\-ment fortement rigides. On calcule
le deuxi\`{e}me groupe
de cohomologie $\mathbf{H}^2_{CE}(\mathfrak{g},\Sg)$ \`{a} l'aide du
th\'{e}or\`{e}me de Hochschild-Serre (\ref{TheoH.S}) appliqu\'{e} \`{a} tous les
groupes
$\mathbf{H}_{CE}^2(\mathfrak{g},\mathcal{S}^k\mathfrak{g})$.

\subsection{L'alg\`{e}bre de Lie de dimension un}

\noindent Evidemment, $\mathbf{H}_{CE}^2(\korps,\mathcal{S}^k\korps)=\{0\}$,
 alors
\bprop \lbl{PUnForRig}
 L'alg\`{e}bre de Lie de dimension $1$ est infinit\'{e}simalement fortement rigide.
\eprop

\subsection{Alg\`{e}bres de Lie semi-simples}

D'apr\`{e}s les lemmes de Whitehead, le premier et le deuxi\`{e}me groupe de
cohomologie d'une alg\`{e}bre de
Lie semi-simple $\mathfrak{g}$  \`{a} valeurs dans tout $\korps$-module de
dimension finie sont triviaux. Par
cons\'{e}quent, $H_{CE}^{2}(\mathfrak{g},\mathcal{S}^{k}\mathfrak{g})=\{0\}$,
$\forall k\in\nat$. Ceci assure, pour le cas particulier $k=1$,
la rigidit\'{e} de l'alg\`{e}bre de Lie. D'apr\`{e}s ce qui pr\'{e}c\`{e}de
on a le th\'{e}or\`{e}me classique (voir \cite[p.~426]{Kas}):
\bsat \lbl{PSemiSForRig}
 Toute alg\`{e}bre de Lie semisimple de dimension finie est infinit\'{e}simalement
 fortement rigide.
\esat

\subsection{Sommes directes de semi-simples et d'autres}

\bsat \lbl{TSemiSFForRig}
Soit $\mathfrak{f}$ une alg\`{e}bre de Lie de dimension finie qui est
infinit\'{e}simalement fortement rigide
%dont
%$\mathbf{H}^2_{CE}(\mathfrak{f}, \mathcal{S}^k\mathfrak{f}) = \{0\}$ quel que
%soit $k\in\nat$
et $\mathfrak{s}$ une alg\`{e}bre de Lie semi-simple de dimension finie.\\
Alors la somme directe $\mathfrak{g}:=\mathfrak{f}\oplus \mathfrak{s}$ est
infinit\'{e}simalement fortement rigide.
\esat
\bbew
  D'apr\`{e}s le th\'{e}or\`{e}me de Hochschild-Serre (\ref{TheoH.S}) et les lemmes de
  Whitehead
 on a
 \[
  \mathbf{H}^2_{CE}(\mathfrak{g},\Sg)
         \cong \mathbf{H}^2_{CE}(\mathfrak{f},\Sg_{a})^\mathfrak{s}
         \cong \mathbf{H}^2_{CE}(\mathfrak{f},\mathcal{S}\mathfrak{f}\otimes
                                       \mathcal{S}\mathfrak{s})^\mathfrak{s}
         \cong \mathbf{H}^2_{CE}(\mathfrak{f},\mathcal{S}\mathfrak{f}_{a})\otimes
                                       (\mathcal{S}\mathfrak{s})^\mathfrak{s}
         \cong\{0\}
 \]
 et $\mathfrak{g}$ est infinit\'{e}simalement fortement rigide.
\ebew

\noindent En particulier:
\bcor\lbl{Prop3.1}
La somme directe d'une alg\`{e}bre
de Lie semi-simple $\mathfrak{s}$ de dimension finie et de l'alg\`{e}bre de Lie
de dimension $1$ est infinit\'{e}simalement fortement rigide.
\ecor

\noindent Par exemple, l'alg\`{e}bre de Lie $\mathfrak{gl}(m,\korps)$ est
infinit\'{e}simalement fortement rigide quel que soit $m\in\nat$.

\subsection{Alg\`{e}bre de Lie non ab\'{e}lienne de dimension $2$}

\noindent On d\'{e}signe par $\mathfrak{r}_{2}$ l'alg\`{e}bre de Lie r\'{e}soluble
engendr\'{e}e par $X,Y$ tel que
\[
      [X,Y]=Y.
\]
Celle-ci est l'alg\`{e}bre de Lie du groupe affine de la droite.

\bprop
L'alg\`{e}bre de Lie $\mathfrak{r}_{2}$ est infinit\'{e}simalement fortement rigide.
%de Lie du groupe affine de la droite. Alors:
%\[
%    \mathbf{H}_{H}^{2}(\mathcal{U}\mathfrak{r}_{2},
%         \mathcal{U}\mathfrak{r}_{2})
%                \simeq
%    \mathbf{H}_{CE}^{2}(\mathfrak{r}_{2},
%         \mathcal{S}\mathfrak{r}_{2})=\{0\}
%\]
%\noindent L'alg\`{e}bre de Lie $\mathfrak{r}_{2}$ est donc fortement rigide.
\eprop
\bbew
%D'apr\`{e}s le th\'{e}or\`{e}me de Cartan-Eilenberg et l'isomorphie des
%$\mathfrak{r}_2$-modules $\mathcal{U}\mathfrak{r}_{2}$ et
%$\mathcal{S}\mathfrak{r}_{2}$ il vient
%\[
%   \mathbf{H}_{H}^{2}(\mathcal{U}\mathfrak{r}_{2},
%         \mathcal{U}\mathfrak{r}_{2})
%                \simeq
%    \mathbf{H}_{CE}^{2}(\mathfrak{r}_{2},
%         \mathcal{U}\mathfrak{r}_{2})
%                 \simeq
%    \mathbf{H}_{CE}^{2}(\mathfrak{r}_{2},
%         \mathcal{S}\mathfrak{r}_{2}).
%\]
Posons, $\mathfrak{r}_{2}=\mathfrak{t}\oplus\mathfrak{n}$ o\`{u}
$\mathfrak{t}:=\korps X$ est une sous-alg\`{e}bre r\'{e}ductive et
$\mathfrak{n}:=\korps Y$ est le nilradical de $\mathfrak{r}_{2}$.
Puisque $\mathfrak{t}$ et $\mathfrak{n}$ sont de dimension $1$, il vient
que
\[
  \mathbf{H}^2_{CE}(\mathfrak{n},\mathcal{S}\mathfrak{r}_{2})=\{0\},~~
  \mathbf{H}^2_{CE}(\mathfrak{t},\korps)=\{0\},~~~
  \mathbf{H}^1_{CE}(\mathfrak{t},\korps)=\korps,
\]
et le th\'{e}or\`{e}me de Hochschild-Serre entra\^{\i}ne que
\beq
  \mathbf{H}^2_{CE}(\mathfrak{r}_2,\mathcal{S}\mathfrak{r}_2)\simeq
     \mathbf{H}^1_{CE}(\mathfrak{n},\mathcal{S}\mathfrak{r}_2)^{\mathfrak{t}}.
\end{equation}
Il est clair que l'espace des $1$-cocycles
$\mathbf{Z}^1_{CE}(\mathfrak{n},\mathcal{S}\mathfrak{r}_2)$ est \'{e}gal \`{a}
$\mathbf{C}^1_{CE}(\mathfrak{n},\mathcal{S}\mathfrak{r}_2)$,  l'espace
des $1$-cocha\^{\i}nes. Soit $F$ une telle $1$-cocha\^{\i}ne et soit $f:=F(Y)\in
\mathcal{S}\mathfrak{r}_2$. $F$ est $\mathfrak{t}$-invariante si
$0=X.\big(F(Y)\big)-F([X,Y])=X.f-f$. Il existe un unique polyn\^{o}me
$\phi\in\korps[X,Y]$ tel que $f=\phi(X,Y)$. Puisque
$X.f=Y\frac{\partial \phi(X,Y)}{\partial Y}$, l'invariance de $F$ entra\^{\i}ne
que $\phi(X,Y)=\chi(X)Y$ pour un polyn\^{o}me $\chi\in\korps[X]$.

D'un autre cot\'{e}, soit $G=\psi(X,Y)\in \mathcal{S}\mathfrak{r}_2
=\mathbf{C}^0_{CE}(\mathfrak{n},\mathcal{S}\mathfrak{r}_2)$ o\`{u}
$\psi\in\korps[X,Y]$. $G$ est invariant par $\mathfrak{t}$ si et seulement si
$0=X.G=Y\frac{\partial \psi(X,Y)}{\partial Y}$, c'est \`{a} dire, si et seulement si $\psi(X,Y)=
\eta(X)$ pour un polyn\^{o}me $\eta\in\korps[X]$. Alors $(\delta_{CE}G)(Y)
=Y.G=-Y\frac{\partial \eta(X)}{\partial X}$.

Il est \'{e}vident que pour chaque polyn\^{o}me $\chi\in\korps[X]$, il existe une
primitive $\eta\in\korps[X]$ telle que
$-\frac{\partial \eta(X)}{\partial X}=\chi(X)$. Il s'ensuit que le groupe
de cohomologie
$\mathbf{H}^1_{CE}(\mathfrak{n},\mathcal{S}\mathfrak{r}_2)^{\mathfrak{t}}$
s'annule, donc $\mathfrak{r}_2$ est infinit\'{e}simalement fortement rigide.
\ebew

\noindent A l'aide du th\'{e}or\`{e}me \ref{TSemiSFForRig} on obtient :
\bcor
 La somme directe de $\mathfrak{r}_2$ et d'une alg\`{e}bre de Lie semisimple
 de dimension finie $\mathfrak{s}$ est infinit\'{e}simalement fortement rigide.
\ecor

\subsection{L'alg\`{e}bre $\mathfrak{gl}(m,\korps)\times \korps^m$}

Soit $m$ un entier strictement positif,
$\mathfrak{a}$ l'alg\`{e}bre de Lie $\mathfrak{gl}(m,\korps)$,
$\mathfrak{n}$ l'espace vectoriel $\korps^m$ et
$\mathfrak{g}$ la somme semi-directe $\mathfrak{a}\oplus \mathfrak{n}$
(o\`{u} le crochet est d\'{e}fini par
$[\phi + v,\phi'+ v']:= [\phi,\phi']+\phi v'-\phi' v$). On a le th\'{e}or\`{e}me
suivant qui g\'{e}n\'{e}ralise l'exemple pr\'{e}c\'{e}dent
$\mathfrak{r}_2=\mathfrak{gl}(1,\korps)\times \korps^1)$:

\bsat\lbl{TInfFortRigGlmKm}
 Avec les d\'{e}finitions mentionn\'{e}es ci-dessus on a
 \[
  \mathbf{H}^0_{CE}(\mathfrak{n},\Sg)^\mathfrak{a}
             \cong \korps
              \mathrm{~~~et~~~}
    \mathbf{H}^k_{CE}(\mathfrak{n},\Sg)^\mathfrak{a}\cong \{0\}
    ~~~~\forall k\in\nat, k\geq 1.
 \]
 Il s'ensuit que
 \[
    \mathbf{H}^k_{CE}(\mathfrak{g},\Sg)\cong
     \mathbf{H}^k_{CE}(\mathfrak{a},\korps).
 \]
 En particulier, $\mathbf{H}^2_{CE}(\mathfrak{g},\Sg)\cong \{0\}$, donc
 $\mathfrak{g}$ est infinit\'{e}simalement fortement rigide.
\esat
\bbew
 Puisque $\mathfrak{a}$ est une sous-alg\`{e}bre r\'{e}ductive dans
 $\mathfrak{g}$, le th\'{e}or\`{e}me de Hochschild-Serre \ref{TheoH.S} s'applique
 et donne la d\'{e}composition
 \[
    \mathbf{H}^k_{CE}(\mathfrak{g},\Sg)
     \cong
       \bigoplus_{r=0}^k
       \mathbf{H}^{k-r}_{CE}(\mathfrak{a},\korps)\otimes_\korps
       \mathbf{H}^{r}_{CE}(\mathfrak{n},\Sg)^\mathfrak{a}.
 \]
 Pour calculer les groupes de cohomologie
 $\mathbf{H}^{r}_{CE}(\mathfrak{n},\Sg)^\mathfrak{a}$ il faut consid\'{e}rer
 les espaces des $r$-cocha\^{\i}nes $\mathfrak{a}$-invariantes,
 $\mathbf{Hom}_\korps(\Lambda^r\mathfrak{n}, \Sg)^\mathfrak{a}$.
 Evidemment, $\Sg\cong \mathcal{S}\mathfrak{a}\otimes_\korps
  \mathcal{S}\mathfrak{n}$, et on a
$\mathbf{Hom}_\korps(\Lambda^r\mathfrak{n},\mathcal{S}^l\mathfrak{g})
     = \bigoplus_{i=0}^l
     \mathbf{Hom}_\korps(\Lambda^r\mathfrak{n},\mathcal{S}^i\mathfrak{a}
                         \otimes_\korps\mathcal{S}^{l-i}\mathfrak{n})$.
 Puisque $\mathfrak{a}=\mathbf{Hom}_\korps(\mathfrak{n},\mathfrak{n})$,
 alors $\mathcal{S}^i\mathfrak{a}\subset
 \mathbf{Hom}_\korps(\mathfrak{n},\mathfrak{n})^{\otimes i}
 = \mathbf{Hom}_\korps
 (\mathfrak{n}^{\otimes i},\mathfrak{n}^{\otimes i})$,
 donc
 \[
   \mathbf{Hom}_\korps(\Lambda^r\mathfrak{n},
     \mathcal{S}^i\mathfrak{a}\otimes_\korps
       \mathcal{S}^{l-i}\mathfrak{n})^\mathfrak{a}
      \subset
      \mathbf{Hom}_\korps(\mathfrak{n}^{\otimes (i+r)},
       \mathfrak{n}^{\otimes (i+l-i)})^\mathfrak{a}.
 \]
 Soit $f\in\mathbf{Hom}_\korps(\mathfrak{n}^{\otimes (i+r)},
       \mathfrak{n}^{\otimes l})^\mathfrak{a}\cong
 \big(\mathfrak{n}^{\otimes l}\otimes_{\korps}
 \mathfrak{n}^{*\otimes (i+r)}\big)^\mathfrak{a}$. En particulier, $f$ est
       invariant par l'application identit\'{e} $\mathbf{1}\in\mathfrak{a}$.
 Puisque $\mathbf{1}.g=lg$ pour tout $g\in \mathfrak{n}^{\otimes l}$ et
 $\mathbf{1}.h=-(r+i)h$ pour tout $h\in \mathfrak{n}^{*\otimes (r+i)}$,
 alors $f$ ne peut \^{e}tre non nul que si $r+i=l$, donc
 $\mathbf{Hom}_\korps(\Lambda^r\mathfrak{n},
 \mathcal{S}^l\mathfrak{g})^\mathfrak{a}
 = \mathbf{Hom}_\korps(\Lambda^r\mathfrak{n},\mathcal{S}^{l-r}\mathfrak{a}
     \otimes \mathcal{S}^{r}\mathfrak{n})^\mathfrak{a}$. Puisque
 l'alg\`{e}bre commutative $\Sg$ s'identifie de fa\c{c}on naturelle \`{a} l'alg\`{e}bre
 des polyn\^{o}mes sur l'espace dual $\mathfrak{a}^*\times \mathfrak{n}^*$, alors
 l'application
 lin\'{e}aire $\Phi^l_r$ suivante est une surjection $\mathfrak{a}$-\'{e}quivariante:
 \beas
  \Phi^l_r:\mathbf{Hom}_\korps(\mathfrak{n}^{\otimes l},
       \mathfrak{n}^{\otimes l}) & \ra &
    \mathbf{Hom}_\korps(\Lambda^r\mathfrak{n},\mathcal{S}^{l-r}\mathfrak{a}
     \otimes \mathcal{S}^{r}\mathfrak{n}): \\
     \psi_1\otimes\cdots \otimes\psi_l & \mapsto &
     \Big( (\phi,\alpha)\mapsto \big(
            v_1\wedge\cdots\wedge v_r\mapsto
     \mathrm{trace}(\psi_1\phi)\cdots\mathrm{trace}(\psi_{l-r}\phi) \\
     & & ~~~~~~~~~~~~~~~
%     {\textstyle\frac{1}{r!}}
     \langle
     (\alpha\psi_{l-r+1})\wedge\cdots\wedge(\alpha\psi_l),
       v_1\otimes\cdots\otimes v_r
       \rangle \big)\Big)
 \eeas
 o\`{u} on a utilis\'{e} $\mathfrak{a}^*\cong \mathfrak{a}$, et
 $\phi,\psi_1,\ldots,\psi_l\in\mathfrak{a}$, $\alpha\in\mathfrak{n}^*$ et
 $v_1,\ldots,v_r\in\mathfrak{n}$. Gr\^{a}ce \`{a} l'\'{e}quivariance de $\Phi^l_r$ il
 vient que l'espace $\Phi^l_r\big(\mathbf{Hom}_\korps(\mathfrak{n}^{\otimes l},
       \mathfrak{n}^{\otimes l})^\mathfrak{a}\big)$ est \'{e}gal \`{a} l'espace
 $\mathbf{Hom}_\korps(\Lambda^r\mathfrak{n},\mathcal{S}^{l-r}\mathfrak{a}
     \otimes \mathcal{S}^{r}\mathfrak{n})^\mathfrak{a}$.

 \noindent L'espace des $\mathfrak{a}$-invariants,
 $\mathbf{Hom}_\korps(\mathfrak{n}^{\otimes l},
 \mathfrak{n}^{\otimes l})^\mathfrak{a}$ est d\'{e}termin\'{e} par le
 th\'{e}or\`{e}me classique de H.~Weyl (voir par exemple \cite[p.29]{DC71}
 pour le cas d'un corps alg\'{e}briquement clos ou \cite[p.214]{KMS93}
 pour le cas d'un corps ordonn\'{e} suivant un raisonnement de G.B.Gurevich,
 1948): soit $S_l$ le groupe sym\'{e}trique d'ordre $l$, et pour une
 permutation $\sigma\in S_l$ soit $\hat{\sigma}$ l'\'{e}l\'{e}ment de
 $\mathbf{Hom}_\korps(\mathfrak{n}^{\otimes l},
 \mathfrak{n}^{\otimes l})$ d\'{e}fini par
 \[
    \hat{\sigma}(w_1\otimes\cdots\otimes w_l):= w_{\sigma(1)}\otimes\cdots
        \otimes w_{\sigma(l)}.
 \]
 Alors d'apr\`{e}s le th\'{e}or\`{e}me de H.~Weyl
 \[
   \mathbf{Hom}_\korps(\mathfrak{n}^{\otimes l},
 \mathfrak{n}^{\otimes l})^\mathfrak{a}=
   \{{\textstyle \sum_{\sigma\in S_l}}
          a_\sigma \hat{\sigma}~|~a_\sigma\in\korps\}.
 \]
 Soit $e_1,\ldots,e_m$ une base de $\mathfrak{n}$ et $e^1,\ldots,e^m$ la base
 duale. Alors les \'{e}l\'{e}ments $E^i_j\in\mathfrak{a}$ d\'{e}finis par
 $E^i_j(e_k):=\delta^i_k e_j$ forment une base de $\mathfrak{a}$. Soit
 $\sigma\in S_l$, alors on voit sans peine que
 $\hat{\sigma}=\sum_{i_1,\ldots,i_l=1}^m E^{i_1}_{i_{\sigma(1)}}\otimes
 \cdots\otimes E^{i_l}_{i_{\sigma(l)}}$, alors pour
 $\phi=\sum_{i,j=1}^m \phi^i_j E^j_i \in\mathfrak{a}$
 et $\alpha=\sum_{i=1}^m\alpha_i e^i \in\mathfrak{n}^*$ on a
 \beas
   \Phi^l_r(\hat{\sigma})(\phi,\alpha)
             & = & \sum_{i_1,\ldots,i_l=1}^m
                  \mathrm{trace}(\phi E^{i_1}_{i_{\sigma(1)}})\cdots
                  \mathrm{trace}(\phi E^{i_{l-r}}_{i_{\sigma(l-r)}}) \\
             &   & ~~~~~~(\alpha E^{i_{l-r+1}}_{i_{\sigma(l-r+1)}})
                  \wedge\cdots\wedge
                  (\alpha E^{i_{l}}_{i_{\sigma(l)}}) \\
             & = & \sum_{i_1,\ldots,i_l=1}^m
                      \phi^{i_1}_{i_{\sigma(1)}}\cdots
                    \phi^{i_{l-r}}_{i_{\sigma(l-r)}}
                    ~\alpha_{i_{\sigma(l-r+1)}}\cdots \alpha_{i_{\sigma(l)}}
                         e^{i_{l-r+1}}\wedge\cdots\wedge e^{i_{l}}
 \eeas
 Quand on regarde la d\'{e}composition en cycles de $\sigma$, alors
 on trouve des
 entiers positifs $1\leq a_1\leq \cdots \leq a_p$ et $1\leq b_1 \leq\cdots
 \leq b_r$ tels que $l= a_1+\cdots+a_p+b_1+\cdots+b_r$ et
 \beq \lbl{EqFormCasimir}
   \Phi^l_r(\hat{\sigma})(\phi,\alpha)
    = \epsilon~\mathrm{trace}(\phi^{a_1})\cdots\mathrm{trace}(\phi^{a_p})
        (\alpha\phi^{b_1-1})\wedge\cdots\wedge(\alpha\phi^{b_r-1})
 \end{equation}
 o\`{u} $\epsilon\in\{-1,1\}$. Il est bien connu que les polyn\^{o}mes Casimir
 $\xi^1,\ldots,\xi^m\in \mathcal{S}\mathfrak{a}$ d\'{e}finis par
 $\xi^a:\phi\mapsto \mathrm{trace}(\phi^a)$ engendrent librement l'anneau
 commutatif des invariants $(\mathcal{S}\mathfrak{a})^\mathfrak{a}$,
 voir par exemple \cite{Che55}: c.-\`{a}-d.
 $(\mathcal{S}\mathfrak{a})^\mathfrak{a}$
 est isomorphe \`{a} l'anneau des polyn\^{o}mes en $m$ variables
 $\korps[\xi^1,\ldots,\xi^m]$. Remarquons que $\xi^a$ est de la forme
 $\xi^a=\Phi^l_r\big(\sum_{i_1,\ldots,i_a=1}^m
 E^{i_1}_{i_2}\otimes\cdots\otimes E^{i_{a-1}}_{i_a}\otimes
 E^{i_a}_{i_1}\big)$.

 \noindent Calculons $\delta \xi^a$ o\`{u} $\delta$ d\'{e}signe l'op\'{e}rateur cobord
 de Chevalley-Eilenberg pour les cocha\^{\i}nes de $\mathfrak{n}$: on a
 \beas
 (\delta \xi^a)(v)(\phi,\alpha)
     & = &
     -(v.\xi^a)(\phi,\alpha)
     = -\sum_{s=1}^a\mathrm{trace}(\phi E^{i_1}_{i_2})\cdots
                 \langle \alpha, E^{i_s}_{i_{s+1}}v \rangle\cdots
                 \mathrm{trace}(\phi E^{i_a}_{i_1}) \\
     & = & -\sum_{s=1}^a \phi^{i_1}_{i_2}\cdots\phi^{i_{s-1}}_{i_s}
                 \langle \alpha, E^{i_s}_{i_{s+1}}v \rangle
                 \phi^{i_{s+1}}_{i_{s+2}}\cdots
                  \phi^{i_a}_{i_1}
                  = -a\langle \alpha,\phi^{a-1}v\rangle
 \eeas
 Puisque $\delta$ pr\'{e}serve les cocha\^{\i}nes $\mathfrak{a}$-invariantes, alors
 $\delta \xi^a\in \mathbf{Hom}_\korps(\Lambda^1\mathfrak{n},
 \mathcal{S}^{a}\mathfrak{g})^\mathfrak{a}$. D'apr\`{e}s l'\'{e}quation
 (\ref{EqFormCasimir}) il s'ensuit que tout \'{e}l\'{e}ment de
 $F\in \mathbf{Hom}_\korps(\Lambda^r\mathfrak{n},
 \mathcal{S}^{l}\mathfrak{g})^\mathfrak{a}$ est de la forme
 \beq \lbl{EqPuresCasimir}
      F= \sum_{a_1,\ldots,a_r=1}^m
        F_{a_1\cdots a_r}(\xi^1,\ldots,\xi^m)
        \delta \xi^{a_1}\wedge\cdots\wedge \delta \xi^{a_r}
 \end{equation}
 o\`{u} $F^{a_1\cdots a_r}(\xi^1,\ldots,\xi^m)\in\korps[\xi^1,\ldots,\xi^m]$.
 D'un autre cot\'{e}, il est clair que tout \'{e}l\'{e}ment de la forme
 (\ref{EqPuresCasimir}) est une cocha\^{\i}ne $\mathfrak{a}$-invariante
 de $\mathbf{Hom}_\korps(\Lambda^r\mathfrak{n},
 \mathcal{S}^{l}\mathfrak{g})$. Les $\delta \xi^a$ ($1\leq a\leq m$)
 sont ind\'{e}pendants sur l'anneau $\korps[\xi^1,\ldots,\xi^m]$ car une
 \'{e}quation $g_1\delta\xi^1+\cdots+g_m\delta\xi^m=0$ avec $g_1,\ldots,g_m
 \in \korps[\xi^1,\ldots,\xi^m]$, \'{e}valu\'{e}e sur $(\phi,\alpha,v)\in
 \mathfrak{a}\times \mathfrak{n}^*\times \mathfrak{n}$, impliquerait
 $g_1[\phi]~\mathbf{1}+\cdots+g_m[\phi]~\phi^{m-1}=0$ (o\`{u}
 $g_a[\phi]:=g_a\big(\mathrm{trace}(\phi),\ldots,\mathrm{trace}(\phi^m)\big)$),
 donc $g_a[\phi]=0$ pour la partie dense (au sens de Zariski)
 de $\mathfrak{a}$ de tous les
 $\phi\in\mathfrak{a}$ ayant des valeurs
 propres deux-\`{a}-deux diff\'{e}rentes, alors $g_a=0$ quel que soit $1\leq a\leq
 m$.
 Puisque
 \beas
    \delta F & = & \sum_{a_1,\ldots,a_r=1}^m
        \delta \big(F_{a_1\cdots a_r}(\xi^1,\ldots,\xi^m)\big)
        \wedge \delta \xi^{a_1}\wedge\cdots\wedge \delta \xi^{a_r} \\
             & = & \sum_{a_1,\ldots,a_r,a_{r+1}=1}^m
              \frac{\partial F_{a_1\cdots a_r}}{\partial \xi^{a_{r+1}}}
                     (\xi^1,\ldots,\xi^m)
        \delta \xi^{a_{r+1}}\wedge \delta \xi^{a_1}\wedge\cdots\wedge
               \delta \xi^{a_r}
 \eeas
 on voit sans peine que le cocomplexe
 $\big(\mathbf{Hom}_\korps(\Lambda\mathfrak{n},\Sg)^\mathfrak{a},\delta)$
 est isomorphe au cocomplexe de Koszul
 $\big(\korps[\xi^1,\ldots,\xi^m]\otimes \Lambda \mathfrak{n}^*,
 \delta_K\big)$ qui est acyclique (lemme de Poincar\'{e} alg\'{e}brique),
 voir par exemple \cite[p.~141]{Fed96}. Alors la cohomologie
 $\mathbf{H}_{CE}(\mathfrak{n},\Sg)^\mathfrak{a}$ est isomorphe \`{a}
 $\korps$, ce qui montre le premier \'{e}nonc\'{e} du th\'{e}or\`{e}me.

 \noindent D'apr\`{e}s le th\'{e}or\`{e}me de Hochschild-Serre \ref{TheoH.S}, appliqu\'{e} \`{a} la
 sous-alg\`{e}bre simple $\mathfrak{s}:=\mathfrak{sl}(m,\korps)$ de
 $\mathfrak{a}$ on a
 $\mathbf{H}^2_{CE}(\mathfrak{a},\korps)\cong
  \mathbf{H}^2_{CE}(\mathfrak{s},\korps)\oplus
  \mathbf{H}^1_{CE}(\mathfrak{s},\korps) = \{0\}$ gr\^{a}ce aux lemmes de
  Whitehead.
\ebew

\bcor
La somme directe d'une alg\`{e}bre
de Lie semi-simple $\mathfrak{s}$ de dimension finie et de l'alg\`{e}bre de Lie
 $\mathfrak{gl}(m,\korps)\times \korps^m$ est infinit\'{e}simalement fortement rigide.
\ecor

\section{Propri\'{e}t\'{e}s des alg\`{e}bres de Lie fortement rigides}

\subsection{Rigidit\'{e} de l'alg\`{e}bre de Lie}

\bsat \lbl{Theo1.3}
Si $\mathfrak{g}$ est une alg\`{e}bre de Lie, de dimension
finie sur $\korps$, fortement rigide, alors elle est rigide en tant
qu'alg\`{e}bre de Lie.
\esat
\bbew
On suppose que l'alg\`{e}bre enveloppante $\Ug$ de $\mathfrak{g}$ soit rigide,
mais l'alg\`{e}bre de Lie $\mathfrak{g}$
ne soit pas rigide: alors il existerait une d\'{e}formation formelle non triviale
$(\mathfrak{g}[[t]],\mu_{t})$ de $\mathfrak{g}$ avec
$\mu_{t}=\sum_{n=0}^{\infty}\mu_{n}t^{n}$ telle que la classe de
$\mu _{1}$ soit non nulle dans
$\mathbf{H}_{CE}^{2}(\mathfrak{g},\mathfrak{g})$.
Puisque $\mathfrak g$ est de dimension finie, alors le $\korps[[t]]$-module
$\mathfrak{g}[[t]]$ est
isomorphe au module libre $\mathfrak{g}\otimes_{\korps}\korps[[t]]$.
Soit $e_1,\ldots,e_n$ une base de
$\mathfrak{g}$, $y_1,\ldots,y_n$ les images des vecteurs de base dans
$\Ug$ et  $y'_1,\ldots,y'_n$ les images des vecteurs de base dans
$\mathcal{U}\big(\mathfrak{g}[[t]]\big)$. On note $\bullet$
la multiplication dans $\mathcal{U}\big(\mathfrak{g}[[t]]\big)$.
Pour toute suite finie croissante $I=(i_1,\ldots,i_k)$ d'indices dans
$\{1,\ldots,n\}$
posons $y_I:=y_{i_1}\cdots y_{i_k}$ dans $\Ug$ et
$y'_I:=y'_{i_1}\bullet\cdots\bullet y'_{i_k}$ dans
$\mathcal{U}\big(\mathfrak{g}[[t]]\big)$.
D'apr\`{e}s le th\'{e}or\`{e}me de Poincar\'{e}-Birkhoff-Witt (qui s'applique \`{a} la
situation o\`{u} l'alg\`{e}bre de Lie est un module libre sur un anneau
commutatif, (voir \cite[p.271]{Ca-Ei}))
les $y_I$ forment une base de $\Ug$ sur $\korps$ et les
$y'_I$ forment une base de $\mathcal{U}\big(\mathfrak{g}[[t]]\big)$
sur $\korps[[t]]$. On voit que l'application $\Phi:\Ug\otimes_\korps
\korps[[t]]\ra \mathcal{U}\big(\mathfrak{g}[[t]]\big)$ d\'{e}finie par
$\Phi(y_I):=y'_I$ est un isomorphisme de $\korps[[t]]$-modules.
Soit $\pi_t:\Ug\otimes_\korps \korps[[t]]~\times~\Ug\otimes_\korps
\korps[[t]]~\ra \Ug\otimes_\korps \korps[[t]]$ la multiplication sur le module
$\Ug\otimes_\korps\korps[[t]]$ induite par $\bullet$ et $\Phi$, c.-\`{a}-d.
$\pi_t(a,b):= \Phi^{-1}\big(\Phi(a)\bullet\Phi(b)\big)$. La restriction
de $\pi_t$ aux \'{e}l\'{e}ments de $\Ug\times \Ug$ d\'{e}finit une application
$\korps$-bilin\'{e}aire $\Ug\times \Ug\ra \Ug\otimes_\korps \korps[[t]]\subset
\Ug[[t]]$ que l'on note aussi par $\pi_t$, \`{a} savoir
$\pi_t(u,v)=\sum_{n=0}^{\infty}t^{n}\pi_{n}(u,v)$ quels que soient $u,v\in
\Ug$ o\`{u} $\pi_{n}\in\mathbf{Hom}_\korps(\Ug\otimes \Ug,\Ug)$.
L'associativit\'{e} de $\pi_t$ sur trois \'{e}l\'{e}ments $u,v,w\in\Ug$ entra\^{\i}ne
les \'{e}quations $\sum_{s=0}^r\Big(\pi_s\big(\pi_{r-s}(u,v),w\big)-
\pi_s\big(u,\pi_{r-s}(v,w)\big)\Big)=0$ quel que soit $r\in\nat$.
Par cons\'{e}quent, $\pi_t$ d\'{e}finit une multiplication
$\korps[[t]]$-bilin\'{e}aire associative sur le $\korps[[t]]$-module
$\Ug\hspace{0.5mm}[[t]]$
(qui contient $\Ug\otimes_\korps \korps[[t]]$ comme sous-module dense par
rapport \`{a} la topologie $t$-adique) de fa\c{c}on usuelle:
\[
  \pi_t\left(\sum_{s}^\infty t^s u_s,\sum_{s'=0}^\infty t^{s'} v_{s'}\right)
      := \sum_{r=0}^\infty t^r \sum_{\stackrel{s,s',s''\geq 0}{s+s'+s''=r}}
                     \pi_{s''}(u_s,v_{s'})
\]
En particulier, l'application $\pi_0$ d\'{e}finit une multiplication
associative sur l'espace vectoriel $\Ug$, et
$(\Ug\hspace{0.5mm}[[t]],\pi_t)$ est une d\'{e}formation
formelle associative de $(\Ug,\pi_0)$. Pour toutes suites finies croissantes
 $I,J$ on a
$\pi_0(y_I,y_J)=\Phi^{-1}(y'_I\bullet y'_J)|_{t=0}$: en `ordonnant bien'
le produit $y'_I\bullet y'_J$ on obtient des combinaisons lin\'{e}aires des
$y'_K$ o\`{u} la suite finie croissante $K$ est de longueur inf\'{e}rieure ou \'{e}gale
\`{a} la somme des
longueurs de $I$ et de $J$ dont les coefficients ne contiennent que les
constantes de structure du crochet $\mu_0$ de l'alg\`{e}bre de Lie $\mathfrak{g}$
car $t=0$. Donc $\pi_0$ est \'{e}gale \`{a} la multiplication de l'alg\`{e}bre
enveloppante $\Ug$ et $(\Ug\hspace{0.5mm}[[t]],\pi_t)$ est donc une d\'{e}formation
formelle associative de l'alg\`{e}bre enveloppante $\Ug$.
Il s'ensuit que $\pi_{1}$ est un $2$-cocycle de Hochschild de $\Ug$, et
la restriction de $\pi_{1}$ sur $X,Y\in\mathfrak{g}$ v\'{e}rifie
\beq\lbl{Eq(3)}
    \mu_{1}(X,Y)=\pi_{1}(X,Y)-\pi_{1}(Y,X) ~~~ \forall X,Y\in \mathfrak{g}.
\end{equation}
car l'alg\`{e}bre de Lie $(\mathfrak{g}[[t]],\mu_t)$ est une sous-alg\`{e}bre de
Lie de $\mathcal{U}(\mathfrak{g}[[t]])$ qui, d'apr\`{e}s ce qui pr\'{e}c\`{e}de,
peut \^{e}tre consid\'{e}r\'{e}e comme une
sous-alg\`{e}bre associative de $(\Ug\hspace{0.5mm}[[t]],\pi_t)$.
La rigidit\'{e} de $\Ug$ implique qu'il
existe un isomorphisme formel
$\varphi_{t}=\sum_{r=0}^{\infty}\varphi_{r}t^{r}$, avec
$\varphi_{0}=Id_{\Ug}$ et
$\varphi_{n}\in \mathbf{Hom}_\korps(\Ug,\Ug)$ tel que
\[
    \varphi_{t}(\pi_{t}(u,v))=
            \pi_{t}(\varphi_{t}(u),\varphi_{t}(v)),
            ~~~\forall u,v\in \Ug,
\]
ce qui est \'{e}quivalent \`{a}
\beq \lbl{Eq(4)}
    \sum_{r=0}^{\infty}t^{r}\sum_{\stackrel{a,b\geq 0}{a+b=r}}
       \varphi_{a}(\pi_{b}(u,v))
         =\sum_{r=0}^{\infty }t^{r}\sum_{\stackrel{a,b,c\geq0}{a+b+c=n}}
            \pi_{a}(\varphi_{b}(u),\varphi_{c}(v))\quad \quad
                  \forall u,v\in \Ug.
\end{equation}
Pour $r=1$, la relation pr\'{e}c\'{e}dente devient
\beq \lbl{Eq(5)}
   \pi_{1}(u,v)=
     (\delta_{H}\varphi_{1})(u,v)\quad \quad \forall u,v\in\Ug
\end{equation}
o\`{u} $\delta_{H}$ est l'op\'{e}rateur cobord de Hochschild (voir \cite{Ho})
par rapport \`{a} la multiplication $\pi_0$ de l'alg\`{e}bre enveloppante.

Alors les formules (\ref{Eq(3)}) et (\ref{Eq(5)}) impliquent
\bea
   \mu_{1}(X,Y) & = & (\delta_{H}\varphi _{1})(X,Y)
                  -(\delta_{H}\varphi_{1})(Y,X) \nonumber \\
                & = & X\varphi_1(Y)-\varphi(XY)+\varphi(X)Y
                      -Y\varphi_1(X)+\varphi(YX)-\varphi(Y)X \nonumber \\
                & = & (\delta _{CE}\varphi _{1})(X,Y)\quad \quad
                         \forall X,Y\in\mathfrak{g}
\eea
o\`{u} $\delta_{CE}$ est l'op\'{e}rateur cobord de Chevalley-Eilenberg
(voir \cite{Ca-Ei}).

Par cons\'{e}quent $\mu_{1}$ est un cobord de Chevalley-Eilenberg et sa
classe est nulle dans $\mathbf{H}_{CE}^{2}(\mathfrak{g},\mathfrak{g})$,
d'o\`{u} une contradiction.
\ebew

Ce th\'{e}or\`{e}me permet de restreindre la recherche d'alg\`{e}bres
enveloppantes rigides aux seules alg\`{e}bres de Lie rigides.
Pour les propri\'{e}t\'{e}s et la classification des alg\`{e}bres de Lie rigides,
 voir  (\cite{Ca-Di},\cite{Ca1}, \cite{GA}, \cite{AB-Go} et \cite{GR}).
\subsection{Deuxi\`{e}me groupe de cohomologie scalaire}

Dans ce paragraphe, on \'{e}tablit un crit\`{e}re cohomologique qui assure
qu'une alg\`{e}bre de Lie n'est pas fortement rigide. Soit
$\big(\mathfrak{g},[~,~]\big)$ une
alg\`{e}bre de Lie sur l'anneau commutatif $K$.
Soit $\omega \in \mathbf{Z}_{CE}^{2}(\mathfrak{g},K)$
un $2$-cocycle scalaire de $\mathfrak g$. Soit
$\mathfrak{g}_{\omega}=\mathfrak{g}\oplus K c$ une extension centrale
unidimensionnelle de $\mathfrak{g}$ par rapport \`{a} $\omega$ dont le crochet
$[~,~]^\prime$ est d\'{e}fini de la fa\c{c}on usuelle par
\beq \lbl{EqDefExtCent}
   [X+ac,Y+bc]^\prime:=[X,Y]+\omega(X,Y)c~~~\forall X,Y\in\mathfrak{g};
   a,b\in K.
\end{equation}

\bsat \lbl{Theo2.1}
Soit $\mathfrak{g}$ une alg\`{e}bre de Lie de dimension
finie sur le corps $\korps$ telle que son deuxi\`{e}me groupe de cohomologie
scalaire
$\mathbf{H}_{CE}^{2}(\mathfrak{g},\korps)$ soit non trivial. Alors, $\mathfrak{g}$
n'est pas fortement rigide.
\esat
\bbew
 Soit $\omega\in\mathbf{Z}_{CE}^2(\mathfrak{g},\korps)$ un $2$-cocycle
 de classe non nulle et soit $\mathfrak{g}_{t\omega}[[t]]$
 l'extension centrale unidimensionnelle de l'alg\`{e}bre de Lie
 $\mathfrak{g}[[t]]=\mathfrak{g}\otimes_\korps \korps[[t]]$ sur $K=\korps[[t]]$
 par rapport \`{a} $t\omega$ (voir (\ref{EqDefExtCent})). Dans
 l'alg\`{e}bre enveloppante $\mathcal{U}(\mathfrak{g}_{t\omega}[[t]])$ de
 $\mathfrak{g}_{t\omega}[[t]]$ on note la multiplication par $\bullet$ et on
 consid\`{e}re l'id\'{e}al bilat\`{e}re
 $\mathcal{I}:=(1-c')\bullet\mathcal{U}(\mathfrak{g}_{t\omega}[[t]])=
 \mathcal{U}(\mathfrak{g}_{t\omega}[[t]])\bullet (1-c')$ (o\`{u} $c'$ d\'{e}signe
 l'image de $c$ dans $\mathcal{U}\big(\mathfrak{g}_{t\omega}[[t]]\big)$)
 et l'alg\`{e}bre quotient
 $\mathcal{U}_{t\omega}\mathfrak{g}
 :=\mathcal{U}(\mathfrak{g}_{t\omega}[[t]])/\mathcal{I}$.
 %\[
 % \mathcal{U}_{t\omega}\mathfrak{g}
 % =\mathbf{T}(\mathfrak{g})/<X.Y-Y.X-[X,Y]+t\omega (X,Y)1:X,Y\in \mathfrak{g}>
 %\]
Soit $e_1,\ldots,e_n$ une base de
$\mathfrak{g}$ sur $\korps$.  Alors $c,e_1,\ldots,e_n$ est
une base de $\mathfrak{g}_{t\omega}[[t]]$ sur $\korps[[t]]$.
Soient $y_1,\ldots,y_n$ les images des vecteurs de base dans
$\Ug$ et  $c',y'_1,\ldots,y'_n$ les images des vecteurs de base dans
$\mathcal{U}\big(\mathfrak{g}_{t\omega}[[t]]\big)$.
Pour toute suite finie croissante $I=(i_1,\ldots,i_k)$ d'indices dans
$\{1,\ldots,n\}$
posons $y_I:=y_{i_1}\cdots y_{i_k}$ dans $\Ug$ et
$y'_I:=y'_{i_1}\bullet\cdots\bullet y'_{i_k}$ dans
$\mathcal{U}\big(\mathfrak{g}[[t]]\big)$.
D'apr\`{e}s le th\'{e}or\`{e}me de Poincar\'{e}-Birkhoff-Witt (qui s'applique \`{a} la
situation o\`{u} l'alg\`{e}bre de Lie est un module libre sur un anneau
commutatif, (voir \cite[pp.~271]{Ca-Ei}))
les $y_I$ forment une base de $\Ug$ sur $\korps$, et les
$c^{'\bullet i_0}\bullet y'_I$ (o\`{u} $i_0\in\nat$ et $c^{'\bullet i_0}:=1$)
forment une base de $\mathcal{U}\big(\mathfrak{g}_{t\omega}[[t]]\big)$
sur $\korps[[t]]$. En passant \`{a} l'alg\`{e}bre quotient
$\mathcal{U}_{t\omega}\mathfrak{g}$, on voit que
$c^{'\bullet i_0}$ s'identifie \`{a} $1$. En notant la multiplication dans
$\mathcal{U}_{t\omega}\mathfrak{g}$ par $\cdot$ et les images de
$y'_1,\ldots,y'_n$ par la projection canonique par $y''_1,\ldots,y''_n$,
on voit que tout $y'_I$ est envoy\'{e} sur
$y''_I:=y''_{i_1}\cdot\ldots\cdot y''_{i_n}$. Il s'ensuit que tous
les \'{e}l\'{e}ments $y''_I$ forment une base de l'alg\`{e}bre quotient
$\mathcal{U}_{t\omega}\mathfrak{g}$. Comme dans la d\'{e}monstration
du th\'{e}or\`{e}me
\ref{Theo1.3} pr\'{e}c\'{e}dent, l'application $\Phi:\Ug\otimes_\korps \korps[[t]]\ra
\mathcal{U}_{t\omega}\mathfrak{g}$ donn\'{e}e par $y_I\mapsto y''_I$
d\'{e}finit un isomorphisme de $\korps[[t]]$-modules libres. De fa\c{c}on
compl\`{e}tement analogue \`{a} celle de la d\'{e}monstration pr\'{e}c\'{e}dente on montre
que la multiplication induite sur $\Ug\otimes_\korps \korps[[t]]$ par
la multiplication $\cdot$ de $\mathcal{U}_{t\omega}\mathfrak{g}$ et
$\Phi$ d\'{e}finit une suite d'applications
$\pi_{t}=\sum_{r=0}^{\infty }\pi_{r}t^{r}$, o\`{u} $\pi_{r}\in
\mathbf{Hom}_\korps (\Ug\otimes\Ug,\Ug)$ avec les propri\'{e}t\'{e}s suivantes:
(1.) $\pi_t$ d\'{e}finit une d\'{e}formation formelle associative de l'alg\`{e}bre
associative $(\Ug,\pi_0)$,
c.-\`{a}-d. une multiplication $\korps[[t]]$-bilin\'{e}aire sur
le $\korps[[t]]$-module $\Ug\hspace{0.5mm}[[t]]$ (qui contient
le $\korps[[t]]$-module $\Ug\otimes_\korps \korps[[t]]$ comme sous-espace
dense par rapport \`{a} la topologie $t$-adique), et (2.) $\pi_0$ est la
multiplication usuelle de l'alg\`{e}bre enveloppante $\Ug$ de $\mathfrak{g}$.
Par cons\'{e}quent, $\pi_{1}$ est un  $2$-cocycle de Hochschild de $\Ug$, et sur
tous $X,Y\in \mathfrak{g}\subset \Ug$ on a la
relation: $\omega(X,Y)1=\pi_{1}(X,Y)-\pi_{1}(Y,X)$ car
l'alg\`{e}bre de Lie $\mathfrak{g}_{t\omega}[[t]]$ s'injecte dans
l'alg\`{e}bre quotient $\mathcal{U}_{t\omega}\mathfrak{g}$, et donc dans
$\Ug\otimes_{\korps}\korps[[t]]\subset \Ug\hspace{0.5mm}[[t]]$.

Supposons que $\Ug$ soit rigide, donc en particulier, la d\'{e}formation
$\pi_t$ est triviale. Alors il existe un $1$-cocycle de Hochschild
$\varphi_{1}\in \mathbf{C}_{H}^{1}(\Ug,\Ug)$ tel
que $\pi_{1}=\delta_{H}(\varphi_{1})$. Par suite, il vient
$\forall X,Y\in \mathfrak{g}$:
\[
 \omega (X,Y)1=\pi_{1}(X,Y)-\pi_{1}(Y,X)
                  =\delta_{H}(\varphi_{1})(X,Y)-
                   \delta_{H}(\varphi_{1})(Y,X)
                   =\delta_{CE}(\varphi_{1})(X,Y).
\]
On en d\'{e}duit que $\omega$ est un cobord de Chevalley-Eilenberg et sa
classe est nulle dans $\mathbf{H}_{CE}^{2}(\mathfrak{g},\korps)$,
d'o\`{u} la contradiction.
\ebew

\subsection{Quelques contre-exemples}

\smallskip
Soit $\mathfrak{g}$ une alg\`{e}bre de Lie rigide, selon \cite{Ca1}, elle se
d\'{e}compose sous la forme $\mathfrak{g}=\mathfrak{s}\oplus
            \mathfrak{t}\oplus \mathfrak{n}$,
  o\`{u} $\mathfrak{s}$ est une sous-alg\`{e}bre de Levi, $\mathfrak{n}$ est le
  nilradical et $\mathfrak{t}$ est un tore
  (une sous-alg\`{e}bre ab\'{e}lienne r\'{e}ductive) et  telle que le radical de
  $\mathfrak{g}$ soit \'{e}gal \`{a} $\mathfrak{t}\oplus\mathfrak{n}$ comme espace
  vectoriel.\\ On donne ci-apr\`{e}s une liste d'alg\`{e}bres de Lie
satisfaisant le crit\`{e}re pr\'{e}c\'{e}dent.

\bcor\lbl{Cor2.2}
 Les alg\`{e}bres de Lie suivantes ne sont pas fortement rigides:
 \ben
  \item
  %(i)
  Toute alg\`{e}bre de Lie nilpotente $\mathfrak{g}$ de dimension strictement
  sup\'{e}rieure \`{a} un.
  \item
  %(ii)
  Toute alg\`{e}bre de Lie  $\mathfrak{g}=\mathfrak{s}\oplus
            \mathfrak{t}\oplus \mathfrak{n}$ dont le tore $\mathfrak{t}$
  est de dimension sup\'{e}rieure ou \'{e}gale \`{a} $2$. En particulier,
  les sous-alg\`{e}bres de Borel d'une alg\`{e}bre de Lie semi-simple
  non isomorphe \`{a} $\mathfrak{sl}(2,\korps)$ ne sont pas fortement rigides.
 % Le suivant est douteux: gl(n,k)+k^2 est une sousalg\`{e}bre parabolique
 % de sl(n+1,k)
%  \item
%  %(iii)}
%  Toute sous-alg\`{e}bre parabolique $\mathfrak{p}$
%  d'une alg\`{e}bre de Lie semi-simple $\mathfrak{g}$ telle que
%  $\mathfrak{p}$ est diff\'{e}rente de $\mathfrak{g}$ et de
%  l'alg\`{e}bre de Lie non ab\'{e}lienne de dimension $2$.
  \item
  %(iv)
  Les alg\`{e}bres de Lie rigides qui sont des
%  $\mathfrak{g}=\mathfrak{sl}(2,\korps)\oplus \complex^{n}$
  produits semi-directs de $\mathfrak{sl}(2,\korps)$ par
  les $\mathfrak{sl}(2,\korps)$-modules irr\'{e}ductibles $\korps^{n}$
  avec $n$ pair.
 \een
\ecor
\bbew
La premi\`{e}re assertion est une cons\'{e}quence d'un r\'{e}sultat classique de Dixmier
concernant les alg\`{e}bres de Lie nilpotentes
(\cite{Di1}):
\[
   \mathbf{H}_{CE}^{2}(\mathfrak{g},\korps)\neq 0
            \mathrm{~~~ si~} \dim (\mathfrak{g)}\geq 2.
\]
Pour (2)
%et (3),
on applique le th\'{e}or\`{e}me de Hochschild--Serre
et le fait que pour une sous-alg\`{e}bre ab\'{e}lienne on ait
$\mathbf{H}_{CE}^{2}(\mathfrak{t},\korps)\neq \{0\}$. \\
%(4)
(3) Les $\mathfrak{sl}(2,\korps)$-modules irr\'{e}ductibles $\korps^{n}$ avec $n$
pair poss\`{e}dent une forme symplectique $\omega$ invariante par
$\mathfrak{sl}(2,\korps)$, par suite en appliquant le th\'{e}or\`{e}me de
factorisation de Hochschild--Serre on a
\[
  \mathbf{H}_{CE}^{2}(\mathfrak{g},\korps)\simeq
  \mathbf{H}_{CE}^{2}(\korps^{n},\korps)^{\mathfrak{sl}(2,\korps)}
      \neq \{0\}
\]
puisque ce dernier contient $\omega$.
\ebew

\section{Quantification par d\'{e}formation}

Dans ce paragraphe, on rappelle la notion de star-produit associ\'{e} \`{a}
une structure de Poisson, ainsi que le th\'{e}or\`{e}me de Kontsevitsch qui en
assure l'existence. Le star-produit associ\'{e} \`{a} la structure de
Poisson lin\'{e}aire d'une alg\`{e}bre de Lie servira \`{a} d\'{e}former
l'alg\`{e}bre enveloppante, on parle dans ce cas de d\'{e}formation par
quantification.

\subsection{Star-produits} \lbl{SubSec3.1}

Dans la th\'{e}orie de quantification par d\'{e}formation (voir \cite{BFFLS})
on consid\`{e}re des d\'{e}formations formelles particuli\`{e}res d'une alg\`{e}bre associative et
commutative:

%\textbf{D\'{e}finition}
\bdefi
 Soit $\mathcal{A}=\mathcal{C}^{\infty}(M)$
l'alg\`{e}bre associative commutative des fonctions $\mathcal{C}^{\infty}$
\`{a} valeurs r\'{e}elles ou complexes
sur une vari\'{e}t\'{e} $M$ munie d'une structure de Poisson $P$. Le crochet
de Poisson sur $\mathcal A$ qui en r\'{e}sulte est not\'{e} $\{~,~\}$.

Un $*$-produit sur $\mathcal{A}$ est
une d\'{e}formation associative formelle de $\mathcal{A}$
\[
     f*_{\lambda}g=\sum_{r=0}^{\infty}\lambda^r B_{r}(f,g)
\]
telle que:
\ben
 \item
 %$(a)_{*}$\textit{-}
 Le $*$-produit plong\'{e} dans l'espace $\mathcal{A}[[\lambda]]$ des s\'{e}ries formelles
 en $\lambda$ \`{a} coefficients dans $\mathcal{A}$ est associative, c'est \`{a} dire:
 \[
 \forall r\in N,\quad
    \sum_{i=0}^{r}(B_{i}(B_{r-i}(f,g),h)-(B_{i}(f,B_{r-i}(g,h)))=0.
 \]
 \item
 %$(b)$
 $B_{0}(f,g)=fg$ (la multiplication point par point).
 \item
 %$(c)$
 $B_{1}(f,g)-B_1(g,f)=\{f,g\}$ et $B_{r}(f,g)=(-1)^{r}B_{r}(g,f)$.
 \item
 %$(d)$
 $B_{r}(f,1) =B_{r}(1,f)=0$ $\quad \forall r>0$
 \item
 Les $B_{r}$ sont bilin\'{e}aires.
\een
\edefi

La condition $3$ implique que $[f,g]:=
\frac{1}{2\lambda }\left(f*_{\lambda}g-\;g*_{\lambda}f\right)$ est une
d\'{e}\-for\-ma\-tion de la structure de Lie $\{~,~\}$.

\smallskip

\subsection{La formule de Kontsevitch} \lbl{SubSec3.2}

\medskip

Dans le domaine des star-produits, le r\'{e}sultat le plus spectaculaire fut
donn\'{e} par Kontsevitch en 1997, (voir \cite{Ko}):

Soit $M$ une vari\'{e}t\'{e} de Poisson. Soient $(x_{1},\ldots,x_{n})$
des coordonn\'{e}es locales sur $M$. On pose
$\mathcal{A}=\mathcal{C}^{\infty}(M)$.

Soit $P=\frac{1}{2}\sum P^{a,b}\partial_{a}\wedge \partial_{b}$ (avec
$\partial_{i}$ la d\'{e}riv\'{e}e partielle $\frac{\partial }{\partial x_{i}},i=a,b$)
un bivecteur de Poisson sur $M$.

\bsat[M.Kontsevitch 1997]
Il existe un star-produit sur toute vari\'{e}t\'{e} de Poisson $(M,P)$.
\esat

Pour la vari\'{e}t\'{e} de Poisson
$(\real^n,\frac{1}{2}\sum_{a,b=1}^nP^{ab}\partial_a\wedge\partial_b)$
Kontsevitch utilise l'An\-satz suivant pour
l'op\'{e}rateur bidiff\'{e}rentiel $B_r$ du star-produit: soient $f,g\in
\Cinf(\real^n)$, $2r=n_1+\cdots+n_r+M+N$ une partition
de l'entier positif $2r$ en somme d'entiers positifs et soit $\sigma$
une permutation de $\{1,2,\ldots,2r\}$. On note le couple
$\big((n_1,\ldots,n_r,M,N),\sigma\big)$ par $\Gamma_r$, et on d\'{e}finit
l'op\'{e}rateur bidiff\'{e}rentiel par
\bea
  B_{\Gamma_r}(f,g) & := &
      \sum_{a_1,\ldots,a_{2r}=1}^n\left(
       \frac{ \partial^{n_1}P^{a_{\sigma(1)}a_{\sigma(2)}} }
                { \partial x_{a_1}\cdots \partial x_{a_{n_1}} }
                \cdots
       \frac{ \partial^{n_r}P^{a_{\sigma(2r-1)}a_{\sigma(2r)}} }
                { \partial x_{a_{n_1+\cdots+n_{r-1}+1}} \cdots
                       \partial x_{a_{n_1+\cdots+n_r}} }
                \right. \nonumber \\
                    &    &
      \left. \frac{ \partial^{M}f }
              { \partial x_{a_{n_1+\cdots+n_r+1}} \cdots
                \partial x_{a_{n_1+\cdots+n_r+M}} }~~
             \frac{ \partial^{N}g }
              { \partial x_{a_{n_1+\cdots+n_r+M+1}} \cdots
                \partial x_{a_{2r}} } \right). \nonumber \\
                    &    &     \lbl{KontsevitchOp}
\eea
L'op\'{e}rateur $B_r$ s'obtient comme une combinaison lin\'{e}aire particuli\`{e}re
d'op\'{e}\-ra\-teurs du type pr\'{e}c\'{e}dent param\'{e}tr\'{e}e par tous les couples $\Gamma_r$
possibles. Kontsevitch repr\'{e}sente les $\Gamma_r$ par des graphes \`{a} $r+2$
sommets (correspondant \`{a} $r$ structures de Poisson et \`{a} deux fonctions) et
$2r$ arr\^{e}tes (correspondant \`{a} $2r$ d\'{e}riv\'{e}es partielles) dans le demi-plan
sup\'{e}rieur, et les coefficients $w_{\Gamma_r}$ r\'{e}els de la combinaison
lin\'{e}aire s'obtiennent par une int\'{e}gration li\'{e}e \`{a} l'image g\'{e}om\'{e}trique du
graphe, (voir \cite{Ko} ou \cite{AMM99} pour des d\'{e}tails). Dans le cas simple
o\`{u} $\Gamma_r=(0,\ldots,0,M=r,N=r,\sigma)$,
$\sigma(2k-1):=k,\sigma(2k)=r+k$ quel que soit $1\leq k \leq r$ le coeficient
$w_{\Gamma_r}=\frac{1}{r!}$. Donc le
star-produit de Kontsevitch $\star_K$ s'\'{e}crit de la mani\`{e}re suivante:
\beq \lbl{EqDefKontsStar}
   \star_K := \sum_{r=0}^\infty \lambda^r B_r
       :=\sum_{r=0}^\infty\lambda^r
       \sum_{\Gamma_r}w_{\Gamma_r} B_{\Gamma_r}.
\end{equation}

La formule de Kontsevitch s'\'{e}tend facilement aux structures de Poisson
complexes (qui sont des sections dans le fibr\'{e} de bivecteurs complexifi\'{e}),
et
on voit ais\'{e}ment que si $\real^{2n}\cong\complex^n$ et $f,g$ et $P$ sont
holomorphes, alors $f\star g$ sera holomorphe. De plus, si $f,g$ et $P$
sont des applications polyn\^{o}miales complexes, alors il en est de m\^{e}me pour
le star-produit de Kontsevitch $f\star g$. Finalement, si la structure
de Poisson $P$ elle-m\^{e}me est une s\'{e}rie formelle
$P=\sum_{s=0}^\infty t^s P_s$ o\`{u} tous les $P_s$ sont des champs de
bivecteurs tels que $[P,P]=0$, la formule de Kontsevitch donne une
d\'{e}formation de la vari\'{e}t\'{e} de Poisson $(\real^n,P_0)$ \`{a} {\em deux}
param\`{e}tres $\lambda$ et $t$.

\subsection{L'alg\`{e}bre enveloppante $\Ug$ \`{a} la Kontsevitch} \lbl{SubSec3.3}

Soit $(\mathfrak{g},[~,~])$ une alg\`{e}bre de Lie de dimension finie sur
$\mathbb{C}$ et
$\mathfrak{g}^{*}$ son dual alg\'{e}brique. L'alg\`{e}bre sym\'{e}trique
$\Sg$ est identifi\'{e}e \`{a} l'alg\`{e}bre des fonctions polyn\^{o}mes sur
$\mathfrak{g}^{*}$.
Soient $(e_{i})_{i=1...n}$ une base de $\mathfrak{g}$,
et $(e^{i})_{i=1...n}$ sa
base duale et $x =\sum_{i=1}^{n}x_{i}e^{i}\in \mathfrak{g}^{*}$. Soit
$P_0\in \Sg \otimes \Lambda^{2}\mathfrak{g}^{*}$ la structure
de Poisson lin\'{e}aire (\ref{Eq(2)}) provenant du crochet de Lie
de $\mathfrak{g}$.
% , avec
%$P(x)=\frac{1}{2}\sum_{a,b}P^{ab}(x)\partial_{a}\wedge \partial_{b}$,
%avec $\partial_{i}$ la d\'{e}riv\'{e}e partielle $\frac{\partial}{\partial x_{i}}$
%$(i=a,b)$.
%Soit
%\beq
%  \mathsf{E}:=\sum_{i=1}^n x_i\partial_i
%\end{equation}
% le {\em champ d'Euler} sur $\complex^n$. Un star-produit
% $\star=\sum_{r=0}^\infty \lambda^r B_r$ sur $\mathfrak{g}^*$ est dit
% {\em homog\`{e}ne} lorsque
%\beq
%  (\mathsf{E}.B_r)(f,g):=\mathsf{E}\big(B_r(f,g)\big)
%               -B_r(\mathsf{E}f,g)-B_r(f,\mathsf{E}g)=-rB_r(f,g)
%\end{equation}
%quels que soient $f,g\in\mathcal{A}:=\Cinf(\mathfrak{g}^*)$ et $r\in\nat$.

Il y a un star-produit $\star_G$ assez naturel sur $(\mathfrak{g}^*,P_0)$,
d\^{u} \`{a}
Simone Gutt \cite{Gu}: pour tout $\xi\in\mathfrak{g}$ soit
$e_\xi:\mathfrak{g}^*\ra\complex$ la fonction exponentielle
$e_\xi(x):=e^{\langle x,\xi\rangle}$; et pour tous
$\xi,\eta\in\mathfrak{g}$ on d\'{e}finit:
\beq \lbl{EqDefGutt}
    e_\xi\star_G e_\eta := e_{H(\xi,\eta)}
\end{equation}
o\`{u} $H(\xi,\eta)\in\mathfrak{g}[[\lambda]]$ est la s\'{e}rie de
Baker-Campbell-Hausdorff,
\renewcommand{\arraystretch}{0.5}
\bea
   \lefteqn{H(\xi,\eta):= \xi + \eta} \nonumber \\
     & & +\sum_{n=1}^\infty\frac{(-1)^n}{n+1}
          \!\!\!\sum_{\begin{array}{c}
                 \scriptscriptstyle k_1,\ldots,k_n\geq 0 \\
                 \scriptscriptstyle l_1,\ldots,l_n\geq 0 \\
                 \scriptscriptstyle k_i+l_i\geq 1
                \end{array}}
                \!\!\!\!\lambda^{\sum_{i=1}^n k_i+l_i}~~
                \frac{\big(ad(\xi)\big)^{k_1}\big(ad(\eta)\big)^{l_1}\cdots
                      \big(ad(\xi)\big)^{k_n}\big(ad(\eta)\big)^{l_n}}
                      {(k_1+\cdots+k_n+1)k_1!\cdots k_n!l_1!\cdots
                      l_n!}\xi
                    \nonumber  \\
      & &    \lbl{formuleBCH}
\eea
\renewcommand{\arraystretch}{2}
Ce star-produit a les propri\'{e}t\'{e}s suivantes:
\bprop[S.Gutt 1983] \lbl{PGutt}
 Quels que soient les polyn\^{o}mes $f,g$ de $\Sg$ le star-produit
 $f\star_G g$ est un polyn\^{o}me en $\lambda$ et pour $\lambda=1$
 on a
 \[
     f\star_{G1}g:=f\star_G g|_{\lambda=1}
      =\omega^{-1}\left(\omega(f)\omega(g)\right)
 \]
 o\`{u} $\omega$ est l'isomorphisme lin\'{e}aire de sym\'{e}trisation entre $\Sg$ et
 $\Ug$ (\ref{EqDefSymmSgUg}).
 En particulier, l'alg\`{e}bre associative $(\Sg,\star_{G1})$
 est isomorphe \`{a} l'alg\`{e}bre enveloppante $\Ug$.
\eprop
\bbew
 La d\'{e}monstration se d\'{e}duit de la formule (4.2) dans \cite{Gu}, p.256.
\ebew

Il s'av\`{e}re que le star-produit de Gutt est en g\'{e}n\'{e}ral diff\'{e}rent de
celui de Kontsevitch. G.Dito a fait la comparaison:
\bprop[G.Dito 1999] \lbl{PDito}
 Soit $\star_K$ le star-produit de Kontsevitch (\ref{EqDefKontsStar})
 pour la vari\'{e}t\'{e} de Poisson $(\mathfrak{g}^*,P_0)$.\\
 Alors $B_2(f,g)=B_2(g,f)$ quels que soient
 $f,g\in\mathcal{A}$ et
 \beq\lbl{EqCovKonts}
   \xi\star_K \eta - \eta\star_K \xi = \lambda [\xi,\eta]~~~~
     \forall \xi,\eta\in \mathfrak{g}\cong {\mathfrak{g}^*}^*\subset
     \mathcal{A}.
 \end{equation}
 En outre, il existe une \'{e}quivalence $\rho$ entre $\star_G$ et $\star_K$
 (c.-\`{a}-d. $\rho(f\star_G g)=\rho(f)\star_K \rho(g)~~\forall
 f,g\in\mathcal{A}$) qui est de la forme
 \beq \lbl{EqRhoDito}
    \rho=\exp\left(\sum_{r\geq 2}\lambda^r a_r D_r\right)
 \end{equation}
 o\`{u} $a_r\in\complex$ et $D_r$ est un op\'{e}rateur diff\'{e}rentiel d'ordre $r$
 \`{a} coefficients constants qui contient $r$ d\'{e}riv\'{e}es partielles.
\eprop
\bbew
% Soit $r\in\nat$. Puisque $\mathsf{E}$ est une d\'{e}rivation de la
% multiplication point-par-point $B_0$ il vient que $\mathsf{E}.B_0=0$.
% De plus, dans l'op\'{e}rateur $B_{\Gamma_r}$ il y a $r$ structures de Poisson
% $P_0$ et $2r$ d\'{e}riv\'{e}es partielles: puisque $\mathsf{E}P_0^{ij}=P_0^{ij}$
% et $[\mathsf{E},\partial_i]=-\partial_i$ quels que soient $1\leq i,j\leq
% n$ il vient que $\mathsf{E}.B_{\Gamma_r}=-rB_{\Gamma_r}$. La sym\'{e}trie
% de $B_2$ se d\'{e}duit de la formule explicite
% \'{e}quation (3), p.5 de l'article de Dito \cite{Dit99}. Finalement,
% l'\'{e}quation (\ref{EqCovKonts}) est montr\'{e}e dans \cite{Dit99}, Lemma 3,
% p.~9.
Ces \'{e}nonc\'{e}s sont tous d\'{e}montr\'{e}s dans l'article de Dito \cite{Dit99}.
\ebew

%\bbew
% Tous ces \'{e}nonc\'{e}s se trouvent dans \cite{Gu} et \cite{Dit99}. Une formule
% explicite pour $C_2$ (qui est diff\'{e}rent de $B_2$ de Kontsevitch) est
% donn\'{e}e dans \cite{Dit99}, p.8.
%\ebew
\noindent On a le corollaire suivant:
\bcor \lbl{CVersKontUg}
 Soient $f,g\in \Sg$. Alors le star-produit $f\star_K g$ est un polyn\^{o}me
 en $\lambda$ et pour $\lambda = 1$ on obtient donc une multiplication
 associative $\star_{K1}$ sur $\Sg$ d\'{e}finie par
 \beq \lbl{EqDefKontconv}
    f\star_{K1}g := \left(f\star_K g\right)|_{\lambda=1}.
 \end{equation}
 L'alg\`{e}bre associative
 $(\Sg,\star_{K1})$ est isomorphe \`{a} l'alg\`{e}bre enveloppante $\Ug$
 de $\mathfrak{g}$.
\ecor
\bbew
 Soient $\delta_1$ et $\delta_2$ les degr\'{e}s des polyn\^{o}mes $f$ et $g$.
 L'op\'{e}rateur bidiff\'{e}rentiel $B_r$ contient $r$ fois la structure de
 Poisson lin\'{e}aire $P_0$ et $2r$ d\'{e}riv\'{e}es partielles. Donc il restent
 toujours au moins $r$ d\'{e}riv\'{e}es qui sont distribu\'{e}es sur $f$ et $g$.
 Il s'ensuit
 que $B_r(f,g)=0$ si $r> \delta_1+\delta_2$, d'o\`{u} la convergence du
 star-produit $\star_K$ sur les polyn\^{o}mes. On d\'{e}duit de la formule
 (\ref{EqRhoDito}) de Dito que l'op\'{e}rateur diff\'{e}rentiel $\rho_r$
 dans $\rho=\sum_{r=0}^\infty\lambda^r\rho_r$ contient exactement
 $r$ d\'{e}riv\'{e}es partielles: alors $\rho$ converge aussi sur tout polyn\^{o}me
 $f$ et d\'{e}finit ainsi une bijection $\complex$-lin\'{e}aire $\rho_{(1)}:\Sg\ra
 \Sg$ par $\rho_{(1)}(f):=\rho(f)|_{\lambda = 1}$. D'apr\`{e}s la proposition
 \ref{PDito} de Dito $\rho_{(1)}$ d\'{e}finit un isomophisme
 d'alg\`{e}bres associatives $(\Sg,\star_{G1})\ra (\Sg,\star_{K1})$. Alors
 $(\Sg,\star_{K1})$ est isomorphe \`{a} l'alg\`{e}bre enveloppante gr\^{a}ce \`{a} la
 proposition \ref{PGutt} de Gutt.
\ebew

\section{D\'{e}formations des alg\`{e}bres enveloppantes par quantification}
         \lbl{Sec4.}

Dans ce paragraphe on va \'{e}tudier les d\'{e}formations de l'alg\`{e}bre enveloppante
$\Ug$ d'une alg\`{e}bre de Lie complexe $(\mathfrak{g},[~,~])$ qui s'obtiennent
\`{a} l'aide d'une d\'{e}formation polyn\^{o}miale $P_t$ de la structure de Poisson
lin\'{e}aire $P_0$ (\ref{Eq(2)}) sur l'alg\`{e}bre sym\'{e}trique $\Sg$:

\subsection{Le th\'{e}or\`{e}me}

%La structure de $\mathfrak{g}$ permet de d\'{e}finir un crochet de
%Poisson lin\'{e}aire sur $\mathfrak{g}^{*}.$ Soit $(e_{i})_{i=1...n}$ une base
%de $\mathfrak{g,}$ $(e^{i})_{i=1...n}$ sa base duale et $x
%=\sum_{i=1}^{n}x _{i}e^{i}\in \mathfrak{g}^{*}$
%\[
%  \{f,g\}=P_{0}(df,dg)\quad \quad
%     P_{0}=\frac{1}{2}\sum_{i,j}P_{0}^{ij}\partial_{i}\wedge \partial_{j},~~~
%     P^{ij}(x)=\sum_{k}C_{ij}^{k}x_{i}
%\]
%avec $C_{ij}^{k}$ les constantes de structures de $\mathfrak{g}$
%
%On a le crochet de Schouten $\left[ P_{0},P_{0}\right]_{s} =0$
%(voir \cite{Va}). L'alg\`{e}bre $\Sg$ munie du crochet de
%Poisson lin\'{e}aire $P_{0}$, est une alg\`{e}bre de Lie - Poisson {\tt la notion est
%douteuse}.
\bsat \lbl{Theo4.1}
Soient $(\mathfrak{g},[~,~])$ une alg\`{e}bre de Lie de
dimension finie $n$ sur $\complex$, $\mathfrak{g}^{*}$ son dual
alg\'{e}brique, $\Sg$ l'alg\`{e}bre des fonctions
polyn\^{o}miales sur $\mathfrak{g}^{*}$, $\Ug$
l'alg\`{e}bre enveloppante de $\mathfrak{g}$ et $P_{0}$
le champ de bivecteurs de la structure de Poisson lin\'{e}aire
(\ref{Eq(2)}) sur $\mathfrak{g}^{*}$.

On suppose qu'il existe une suite $(P_{a})_{a\in \mathbb{N}}$
de bivecteurs polyn\^{o}miaux
$(P_{a}\in \Sg\otimes \wedge^{2}\mathfrak{g}^{*})$
tels que tous les crochets de Schouten
$\sum_{i+j=a}[P_{i},P_{j}]_s=0$ quels que soient $a\in \mathbb{N}$ et
$P_{1}$ n'est pas cohomologue \`{a} $0$.

Alors $P_{t}=\sum_{a\geq 0}t^{a}P_{a}$ est une
d\'{e}formation non triviale de la structure de Poisson $P_{0}$ et
elle induit une d\'{e}formation non triviale de l'alg\`{e}bre enveloppante
$\Ug$ par la formule de Kontsevitch. En particulier, l'alg\`{e}bre de Lie
$\mathfrak{g}$ n'est pas fortement rigide.
\esat
\bbew
On consid\`{e}re la d\'{e}formation $\hat{\star}_K$
de l'alg\`{e}bre sym\'{e}trique $\Sg[[t]]$ induite par les formules de Kontsevitch
(\ref{KontsevitchOp}) et (\ref{EqDefKontsStar}) pour la structure de
Poisson formelle $P_t$.  Alors pour tout entier positif $r$,
les op\'{e}rateurs bidiff\'{e}rentiels
de Kontsevitch (\ref{KontsevitchOp}) sont des s\'{e}ries formelles en $t$:
\[
 B_{\Gamma_{r}}(f,g)=\sum_{k\geq 0}t^{k}B_{k,\Gamma _{r}}
\]
avec
\bea
  B_{k,\Gamma_r}(f,g) & := &
      \sum_{\stackrel{b_1,\ldots,b_r\geq 0}{b_{1}+...+b_{r}=k}} \nonumber \\
     & &  \sum_{a_1,\ldots,a_{2r}=1}^n
       \left(
       \frac{ \partial^{n_1}P_{b_1}^{a_{\sigma(1)}a_{\sigma(2)}} }
                { \partial x_{a_1}\cdots \partial x_{a_{n_1}} }
                \cdots
       \frac{ \partial^{n_r}P_{b_r}^{a_{\sigma(2r-1)}a_{\sigma(2r)}} }
                { \partial x_{a_{n_1+\cdots+n_{r-1}+1}} \cdots
                       \partial x_{a_{n_1+\cdots+n_r}} }
                \right.   \nonumber \\
                    &    &
      \left. \frac{ \partial^{M}f }
              { \partial x_{a_{n_1+\cdots+n_r+1}} \cdots
                \partial x_{a_{n_1+\cdots+n_r+M}} }~~
             \frac{ \partial^{N}g }
              { \partial x_{a_{n_1+\cdots+n_r+M+1}} \cdots
                \partial x_{a_{2r}} } \right). \nonumber \\
                    &    &     \lbl{KontsevitchOpPlus}
\eea
On pose $B_{r,k}(f,g):=\sum_{\Gamma_{r}}w_{\Gamma_r}B_{k,\Gamma_{r}}(f,g)$,
et le star-produit de $f$ et de $g$ est de la forme suivante:
\[
  f\hat{\star}_K g=\sum_{r,k\geq 0}\lambda ^{r}t^{k}B_{r,k}(f,g),
\]
 et induit donc une d\'{e}formation \`{a} deux param\`{e}tres $\lambda $ et $t$ de
 l'alg\`{e}bre sym\'{e}trique $\Sg$.

On va montrer que cette d\'{e}formation est convergente en $\lambda $ sur
les polyn\^{o}mes: soient $f$ et $g$ deux polyn\^{o}mes dans $\Sg$ tels que
$f$ est homog\`{e}ne de degr\'{e} $\delta_{1}$ et $g$ est homog\`{e}ne de degr\'{e}
$\delta_{2}$. On fixe $r$ et $k$.
Alors le coefficient en $t^{k}$ de $B_{\Gamma _{r}}(f,g)$,
$B_{k,\Gamma _{r}}(f,g)$ consiste en $r-k$
copies de la structure de Poisson lin\'{e}aire $P_{0}$ et $k$ copies des
autres bivecteurs $P_{b_{i}}\in\mathcal{S}\mathfrak{g}\otimes
\Lambda^2\mathfrak{g}^*$ o\`{u} $1\leq b_i\leq k$ quel que soit $b_i$. Alors pour
$k$ fix\'{e} il n'y a qu'un nombre fini $N_k$ des bivecteurs $P_{b_i}$ pour
lesquels $b_i\geq 1$. Soit $d$ le maximum des degr\'{e}s de $P_0$ et de ces
autres bivecteurs.
D'apr\`{e}s la formule (\ref{KontsevitchOpPlus}), le degr\'{e} polyn\^{o}mial $\delta$ de
$B_{k,\Gamma_r}(f,g)$ est major\'{e} de la mani\`{e}re suivante:
\[
 \delta\leq (r-k)1+dk+\delta_1+\delta_2-2r=
            (d-1)k+\delta_1+\delta_2-r
\]
Par cons\'{e}quent, $B_{k,\Gamma_r}(f,g)=0$ si
$r\geq (d-1)k+\delta_{1}+\delta_{2}+1$. Il s'ensuit que cette d\'{e}formation est
convergente en $\lambda$ sur les polyn\^{o}mes $f,g$.
On peut donc fixer $\lambda =1$, et on obtient une multiplication associative
$\hat{\star}_{K1}$ \`{a} un seul param\`{e}tre $t$ sur $\Sg[[t]]$. Pour $k=0$, la
formule (\ref{KontsevitchOpPlus}) ne contient que la structure de Poisson
lin\'{e}aire $P_0$: alors la multiplication $\hat{\star}_{K1}$ est une
d\'{e}formation formelle associative (en $t$) de la multiplication
$\star_{K1}$ (voir (\ref{EqDefKontconv})) sur $\Sg$. Puisque l'alg\`{e}bre
associative complexe $(\Sg,\star_{K1})$ est isomorphe \`{a} l'alg\`{e}bre
enveloppante $\Ug$ d'apr\`{e}s le corollaire \ref{CVersKontUg}, on a donc
construit une
d\'{e}formation formelle associative $\mu_t=\sum_{r=0}^\infty t^r\mu_r$
de l'alg\`{e}bre enveloppante $\Ug$.

Supposons maintenant que $\mu_t$ soit triviale: alors il existerait
une s\'{e}rie $S_t:=1+\sum_{r=1}^\infty t^r S_t$ o\`{u}
$S_r\in\mathbf{Hom}_\complex(\Ug,\Ug)$ telle que $S_t\big(\mu_t(f,g)\big)
=\mu_0\big(S_t(f),S_t(g)\big)$ quels que soient $f,g\in\Ug[[t]]$. En
particulier, pour $r=1$ et $f=\xi,g=\eta$ avec $\xi,\eta\in\mathfrak{g}$
on aurait
\[
   \mu_1(\xi,\eta)=-S_1\big(\mu_0(\xi,\eta\big)+\mu_0\big(S_1(\xi),\eta\big)
                    +\mu_0\big(\xi,S_1(\eta)\big),
\]
alors, gr\^{a}ce \`{a} (\ref{EqCovKonts}):
\beas
  \mu_1^-(\xi,\eta)&:=& \mu_1(\xi,\eta)-\mu_1(\eta,\xi) \\
             & = & -S_1[\xi,\eta]+S_1(\xi)\star_{K1}\eta
                        -\eta\star_{K1} S_1(\xi)
                        -S_1(\eta)\star_{K1}\xi
                        +\xi\star_{K1}S_1(\eta).
\eeas
Puisque la transformation d'\'{e}quivalence de Dito (\ref{EqRhoDito}) en
$\lambda=1$, $\rho_{(1)}$, donne $\rho_{(1)}(\xi)=\xi$ quel que soit
$\xi\in\mathfrak{g}$ on obtient
\beas
  \rho_{(1)}^{-1}\big(\mu_1^-(\xi,\eta)\big) &=&
     -\rho_{(1)}^{-1}\big(S_1[\xi,\eta]\big)
                     +\rho_{(1)}^{-1}\big(S_1(\xi)\big)\star_{G1}\eta
                        -\eta\star_{K1} \rho_{(1)}^{-1}\big(S_1(\xi)\big)
                            \\
                            & &
                        ~~~~~~-\rho_{(1)}^{-1}\big(S_1(\eta)\big)\star_{K1}\xi
                        +\xi\star_{K1}\rho_{(1)}^{-1}\big(S_1(\eta)\big)
                              \\
                              & = &
         -\rho_{(1)}^{-1}\big(S_1[\xi,\eta]\big)
         + [\rho_{(1)}^{-1}\big(S_1(\xi)\big),\eta]
         - [\rho_{(1)}^{-1}\big(S_1(\eta)\big),\xi] \\
         & = &
            \big(\delta_{CE}(\rho_{(1)}^{-1}\circ S_1)\big)(\xi,\eta)
\eeas
car le star-produit de Gutt $\star_{G1}$ est isomorphe \`{a} la multiplication
associative sur $\Sg$ obtenue \`{a} l'aide de l'isomorphisme $\omega$
(\ref{EqDefSymmSgUg}) qui est un morphisme de $\mathfrak{g}$-modules.
Le terme $\mu_1$ est une somme d'op\'{e}rateurs bidiff\'{e}rentiels du type
$B_{k,\Gamma_{r+1}}(\xi,\eta)$ (\ref{KontsevitchOpPlus}) avec $k=1$. Cet
op\'{e}rateur contient $r$ fois la structure de Poisson lin\'{e}aire $P_0$ et
une fois le champs de bivecteurs $P_1$. Si $\epsilon$ est le degr\'{e} maximal
de $\mu_1|_{\mathfrak{g}\times\mathfrak{g}}$, alors le degr\'{e} de
$B_{k,\Gamma_{r+1}}(\xi,\eta)$ est inf\'{e}rieur ou \'{e}gal \`{a} $\epsilon-r$,
donc il est maximal pour $r=0$. Dans ce cas, la composante de degr\'{e}
$\epsilon$ de $B_{k,\Gamma_{r+1}}(\xi,\eta)$ est proportionnelle \`{a}
$P_1(\xi,\eta)$ avec un poids non nul. En consid\'{e}rant la composante
de degr\'{e} $\epsilon$ dans l'\'{e}quation pr\'{e}c\'{e}dente et en notant que
$\rho_{(1)}$ ne change pas la composante du plus haut degr\'{e}, on voit
que $P_1$ est un $2$-cobord de la cohomologie de Chevalley-Eilenberg,
contrairement \`{a} l'hypoth\`{e}se. Alors la d\'{e}formation de l'alg\`{e}bre
enveloppante est non triviale.
\ebew
\begin{remarque}
Les deux th\'{e}or\`{e}mes \ref{Theo1.3} et \ref{Theo2.1} peuvent se d\'{e}duire
\'{e}galement du th\'{e}or\`{e}me \ref{Theo4.1}.
\end{remarque}
\subsection{Applications} \lbl{SubSec4.2}

\medskip On montre que les alg\`{e}bres de Lie rigides $\mathfrak{t}_{1}\oplus
\mathfrak{n}_{5,6}$ et $\mathfrak{sl}(2,\mathbb{C})\oplus \mathfrak{n}_{3}$
ne sont pas fortement rigides en appliquant le th\'{e}or\`{e}me \ref{Theo4.1}.

\subsubsection{Alg\`{e}bre de Lie r\'{e}soluble
$\mathfrak{t}_{1}\oplus \mathfrak{n}_{5,6}$ } \lbl{SubsubSecEx1}

On consid\`{e}re l'alg\`{e}bre de Lie r\'{e}soluble rigide
$\mathfrak{t}_{1}\oplus \mathfrak{n}_{5,6}$ de dimension $6$, d\'{e}finie
relativement \`{a} la base $\{ X_{0},X_{1},X_{2},X_{3},X_{4},X_{5} \}$ par
\bea
     ~ [ X_{0},X_{i} ] & = & iX_{i} \quad i=1,\cdots,5   \\
     ~ [ X_{1},X_{i} ] & = & X_{i+1}\quad i=2,3,4   \\
     ~ [ X_{2},X_{3} ] & = & X_{5}
\eea
les autres crochets \'{e}tant nuls ou d\'{e}duits des pr\'{e}c\'{e}dents par antisym\'{e}trie.

\bprop\lbl{PPoisRes}
%\vspace{0.03in}\textbf{Proposition}
Soit $P\in \Sg\otimes \mathbf{\wedge}^{2}\mathfrak{g}^{*}$
d\'{e}fini par ($\alpha,\beta,\gamma\in\complex^3\setminus\{(0,0,0)\}$):
\bea
    P_{1} & = & \beta X_{2}^{2}
        \frac{\partial}{\partial X_{1}}\wedge \frac{\partial}{\partial X_{3}}
                    \nonumber \\
          &   & +\gamma(-X_{2}X_{3}\frac{\partial}{\partial X_{1}}\wedge
                  \frac{\partial}{\partial X_{4}}
        +X_{2}X_{5}\frac{\partial}{\partial X_{3}}\wedge
                  \frac{\partial }{\partial X_{4}})
        +\alpha X_{1}X_{5}\frac{\partial}{\partial X_{2}}\wedge
                  \frac{\partial}{\partial X_{4}} \nonumber \\
          &   &     \lbl{EqPoisNonlin}
\eea
Alors $[P,P]_{s}=0=[P_{0},P]_{s}$ et $P$ n'est pas cohomologue \`{a} $0$.
En particulier $\mathfrak{g}$ n'est pas fortement rigide.
\eprop
\bbew
Un calcul direct montre que $[P,P]_{s}=0=[P_{0},P]_{s}.$

\noindent Soit $\mathcal{S}_{2}\mathfrak{g}$
l'espace des polyn\^{o}mes quadratiques
(de degr\'{e} 2) de variables
$X_{0}$, $X_{1}$, $X_{2}$, $X_{3}$, $X_{4}$, $X_{5}$. On
d\'{e}termine l'espace des $2$-cochaines quadratiques
$\mathcal{S}^{2}\mathfrak{g}\otimes \mathbf{\wedge}^{2}\mathfrak{g}^{*}$
et l'espace des $1$-cochaines quadratiques
$\mathcal{S}^{2}\mathfrak{g}\otimes \mathfrak{g}^{*}$.
On montre qu'il n'existe aucun \'{e}l\'{e}ment $A$ de
$\mathcal{S}^{2}\mathfrak{g}\otimes \mathfrak{g}^{*}$ tel que
$d_{CE}A=P$. Alors $\mathfrak{g}$ n'est pas fortement rigide d'apr\`{e}s le
th\'{e}or\`{e}me \ref{Theo4.1}.
\ebew

\subsubsection{Produit semi-direct de $\mathfrak{sl}(2,\mathbb{C)}$
par $\mathfrak{n}_{3}$} \lbl{SubsubSecEx2}

On consid\`{e}re l'alg\`{e}bre de Lie parfaite rigide
$\mathfrak{g}=\mathfrak{sl}(2,\mathbb{C)}\oplus \mathfrak{n}_{3}$
produit semi-direct de $\mathfrak{sl}(2,\mathbb{C)}$
par $\mathfrak{n}_{3}$, l'alg\`{e}bre de Heisenberg.

Soient ${h,e,f}$ la base de $\mathfrak{sl}(2,\mathbb{C})$ et ${p,q,z}$
la base de
l'alg\`{e}bre de Heisenberg $\mathfrak{n}_{3}$ v\'{e}rifiant les relations
suivantes :
\[
 \begin{array}{ccclccclcccl}
~[h,e] & = & 2e &, & [h,f]  & = & -2f &, &
                                        [e,f] & = & h  & , \\
~[q,p]  & = & z  &, & [h,q]  & = & q   &, &
                                        [h,p]   & = & -p & , \\
~[e,p]   & = & q  &, & [f,q]   & = & p   &  &
                                              &              &    &
 \end{array}
\]
les autres crochets \'{e}tant nuls ou d\'{e}duits des pr\'{e}c\'{e}dents par antisym\'{e}trie.

L'action induite de $\mathfrak{sl}(2,\mathbb{C})$ sur
$\mathfrak{n}_{3}^{*}$ \`{a} l'aide de la base duale
${p^{*},q^{*},z^{*}}$  de ${p,q,z}$ :
\[
 h\cdot q^{*}=-q^{*},~~ h\cdot p^{*}=p^{*},~~e\cdot q^{*}=-p^{*},~~
 f\cdot p^{*}=-q^{*}
\]
les autres actions \'{e}tant nulles.
En appliquant le th\'{e}or\`{e}me de factorisation de Hochschild-Serre on a :
\[
 \mathbf{H}_{CE}^{2}(\mathfrak{g},\Sg)
\simeq \mathbf{H}_{CE}^{2}(\mathfrak{n}_{3},
        \Sg)^{\mathfrak{sl}(2,\mathbb{C})}
\]

\blem \lbl{Lem1}
Soit $\Phi \in \Sg\otimes \mathfrak{g}^{*}$
(1-cocycle) avec $\Phi =g_{1}q^{*}+g_{2}q^{*}+g_{3}z^{*}$,
et $g_{i}\in \Sg$ on a
\bea
 %$(1)$
 \delta _{CE}(z^{*}) & = & -q^{*}\wedge p^{*} \\
 %$(2)$
 \delta _{CE}(\Phi)  & = & ([q,g_{2}]-[p,g_{1}]-g_{3})
           \cdot q^{*}\wedge p^{*} + [q,g_{3}]\cdot q^{*}\wedge z^{*}
                    +[p,g_{3}]\cdot p^{*}\wedge z^{*} \nonumber \\
                     & &
\eea
\elem

\blem \lbl{Lem2}
Le sous-espace vectoriel $\mathbf{F}$ de $\mathcal{S}^{2}\mathfrak{g}$
d\'{e}fini par
\[
 \mathbf{F}=\{ {g}_{3} \in \mathcal{S}^{2}\mathfrak{g}~:~
              [q,{g}_{3}]=[p,{g}_{3}]=0 \}
\]
est engendr\'{e} par $z^{2}$, $q^{2}-2ez$, $pq+hz$, et $p^{2}+2fz$.
On a $\mathbf{F}\subset \mathcal{S}_{\mathfrak{n}_{3}}$
 avec $\mathcal{S}_{\mathfrak{n}_{3}}$
l'id\'{e}al de $\Sg$ engendr\'{e} par
$\mathfrak{n}_{3}$.
\elem

\bprop
Soit $P\in \Sg\otimes \wedge^{2}
\mathfrak{n}_{3}^{*}$ d\'{e}fini par
\[
     P:=C~q^{*}\wedge p^{*}\mathrm{~~avec~~}
     C:=h^{2}+4ef.
\]
$C$ est l'\'{e}l\'{e}ment de
Casimir quadratique de $\mathcal{S}\mathfrak{sl}(2,\mathbb{C})$.\\
Alors $P$ est $\mathfrak{sl}(2,\mathbb{C})$-invariant,
$[P,P]=0=[P_{0},P]$ et $P$ n'est pas
cohomologue \`{a} $0$ dans
$\mathbf{H}_{CE}^{2}(\mathfrak{n}_{3},
           \Sg)^{\mathfrak{sl}(2,\mathbb{C})}$.
En particulier, $\mathfrak{g}$ n'est pas fortement rigide.
\eprop
\bbew
%\textbf{Preuve }
On v\'{e}rifie facilement que $P$ est $\mathfrak{sl}(2,\mathbb{C})$-invariant,
et $[P,P]=0=[P_{0},P]$.
On suppose qu'il existe un $\Phi =g_{1}q^{*}+g_{2}q^{*}+g_{3}z^{*}\in
\mathcal{S}^{2}\mathfrak{g}\otimes \mathfrak{g}^{*}$ tel que
\[
  d_{CE}(\Phi)  =
    ([q,g_{2}]-[p,g_{1}]-g_{3})~ q^{*}\wedge p^{*}
    +  [q,g_{3}]~q^{*}\wedge z^{*}
    +  [p,g_{3}]~p^{*}\wedge z^{*}
                =  C~q^{*}\wedge p^{*}
\]
d'apr\`{e}s le lemme \ref{Lem1}. Puisque $q^{*}\wedge p^{*}$,
$q^{*}\wedge z^{*}$ ,$q^{*}\wedge p^{*}$ forment une base de
$\wedge^{2}\mathfrak{n}_{3}^{*}$, on a
$g_{3}\in \mathbf{F}\subset \mathcal{S}_{\mathfrak{n}}$
d'apr\`{e}s le lemme \ref{Lem2}.
Or $[q,g_{2}]$, $[p,g_{1}]\in \mathcal{S}_{\mathfrak{n}}$
et par suite $C=([q,g_{2}]-[p,g_{1}]-g_{3})\in \mathcal{S}_{\mathfrak{n}}$,
d'o\`{u} une contradiction.

Alors $P$ n'est pas cohomologue \`{a} $0$ dans
$\mathbf{H}_{CE}^{2}(\mathfrak{n}_{3},\Sg)^{\mathfrak{sl}(2,\mathbb{C})}$ et
$\mathfrak{g}$ n'est pas fortement rigide d'apr\`{e}s le th\'{e}or\`{e}me
\ref{Theo4.1}.
\ebew

\subsection{Classification des alg\`{e}bres fortement rigides de dimension
         inf\'{e}rieure \`{a} six} \lbl{Sec6}

On utilise la classification des alg\`{e}bres de Lie rigides de dimension
inf\'{e}rieure \`{a} sept, donn\'{e}e par Carles et Diakit\'{e} \cite{Ca2} \cite{Ca-Di}
%TCIMACRO{\UNICODE{0xe9} [Ca-Di]}
%BeginExpansion
%T\hskip-.05em\llap{\protect\rule[.7ex]{.6em}{.1ex}}%
%EndExpansion
et les r\'{e}sultats pr\'{e}c\'{e}dents pour classifier les alg\`{e}bres de Lie fortement
rigides.

\bsat \lbl{Theo6.1}
Toute alg\`{e}bre de Lie fortement rigide de dimension $n\leq 6$
est isomorphe \`{a} l'une des alg\`{e}bres de Lie suivantes:
\begin{center}
        $\{0\},~~\korps,~~\mathfrak{r}_{2}, ~~\mathfrak{sl}(2,\korps),~~
           \mathfrak{gl}(2,\korps),~~
            \mathfrak{sl}(2,\korps)\times \mathfrak{r}_{2},~~$\\
        $\mathfrak{sl}(2,\korps)\times \mathfrak{sl}(2,\korps),~~
              \mathfrak{gl}(2,\korps)\times \korps^{2}$.
\end{center}

\esat

\section{Lin\'{e}arisation des structures de Poisson}

Ce probl\`{e}me a \'{e}t\'{e} pos\'{e} par Weinstein \cite[pp. 537]{Wei83} :
Etant donn\'{e} une vari\'{e}t\'{e} de Poisson $(M,P)$ dont la structure de Poisson
$P$ s'annule en un point $x_0\in M$, il
s'agit de voir si il existe un diff\'{e}omorphisme ($\Cinf$ ou formel)
tel que $P$ soit isomorphe \`{a} sa partie d'ordre 1 en $x_0$, $P_0$ (voir \cite{Conn84},
\cite{Conn85}, \cite{Wei87}, \cite{Duf90} et \cite{MZ02} pour des r\'{e}sultats
dans ce domaine). D'un point de vue formel c.\`{a}.d. quand on remplace
$P$ par sa s\'{e}rie de Taylor autour $x_0$ on ne consid\`{e}re que des d\'{e}formations
formelles de $P_0$ qui sont au moins quadratiques. La lin\'{e}arisabilit\'{e}
de l'alg\`{e}bre de Lie $\mathfrak{gl}(m,\korps)\times \korps^m$ qui est fortement rigide,
voir le th\'{e}or\`{e}me \ref{TInfFortRigGlmKm}, a \'{e}t\'{e}
obtenue \`{a} l'aide d'une d\'{e}monstration g\'{e}om\'{e}trique par
Dufour et Zung, voir \cite{DNTZ02}.

\noindent Du theor\`{e}me \ref{Theo4.1}, on d\'{e}duit le r\'{e}sultat suivant:
\bsat
Toute structure de Poisson dont la partie lin\'{e}aire correspond \`{a} une
alg\`{e}bre de  Lie fortement rigide est formellement lin\'{e}arisable.
\esat

\noindent La r\'{e}ciproque de ce th\'{e}or\`{e}me
n'est pas vraie en g\'{e}n\'{e}ral: pour l'alg\`{e}bre de Lie $\mathfrak{g}$ somme semi-directe
$\korps\times \korps^2$, o\`{u} la sous-alg\`{e}bre $\korps$ agit
sur le nilradical ab\'{e}lien $\korps^2$ par des multiples de
l'identit\'{e}, on a que le deuxi\`{e}me groupe de cohomologie
$\mathbf{H}_{CE}^2(\mathfrak{g},\bigoplus_{i=0,i\neq 1}^\infty\mathcal{S}^i\mathfrak{g})$
s'annule, alors elle est lin\'{e}arisable (voir \cite{Duf90}), mais
$\mathfrak{g}$ n'est pas rigide en tant qu'alg\`{e}bre de Lie, car on peut
ais\'{e}ment remplacer l'action de $\korps$ sur $\korps^2$ par une autre
action lin\'{e}aire. Donc $\mathfrak{g}$ n'est pas fortement rigide.

\begin{remarque}
La structure de Poisson $P_{0}+P_{1}$ d\'{e}finie dans (\ref{EqPoisNonlin}) de
la proposition \ref{PPoisRes} n'est pas formellement lin\'{e}arisable parce qu'elle admet
$P_1$ comme d\'{e}formation quadratique non triviale.
\end{remarque}

\end{document}